\documentclass[twocolumn]{autart}   

\usepackage{amsmath,amssymb}
\usepackage{graphicx}
\usepackage{makecell,hhline,tabularx}
\usepackage{algpseudocode}
\usepackage{mathtools,cuted}
\usepackage{color}
\usepackage{tikz}
\usepackage{hyperref}
\usetikzlibrary{arrows.meta}

\newcolumntype{L}[1]{>{\raggedright\let\newline\\\arraybackslash\hspace{0pt}}p{#1}}
\newcolumntype{C}[1]{>{\centering\let\newline\\\arraybackslash\hspace{0pt}}p{#1}}

\DeclareMathSizes{10}{10}{7}{4}

\newcommand{\algrule}[1][.2pt]{\par\vskip.2\baselineskip\hrule height #1\par\vskip.2\baselineskip}

\graphicspath{,{./Figs/}}

\begin{document}

\begin{frontmatter}
	\title{Spectral Bayesian Estimation for General \\ Stochastic Hybrid Systems\thanksref{footnoteinfo}}
	
	\thanks[footnoteinfo]{This paper was not presented at any IFAC meeting. Corresponding author W. Wang.}
	
	\author[SEH]{Weixin Wang}\ead{wwang442@gwu.edu},
	\author[SEH]{Taeyoung Lee}\ead{tylee@gwu.edu}
	
	\address[SEH]{George Washington University, 800 22nd ST NW, Washington, DC, 20052}
	
	\begin{keyword}                           
		hybrid system, stochastic jump processes, filtering technique, state estimation, numerical methods.
	\end{keyword}
    
	\begin{abstract}
		General Stochastic Hybrid Systems (GSHS) have been formulated to represent various types of uncertainties in hybrid dynamical systems.
		In this paper, we propose computational techniques for Bayesian estimation of GSHS.
        In particular, the Fokker-Planck equation that describes the evolution of uncertainty distributions along GSHS is solved by spectral techniques, where an arbitrary form of  probability density of the hybrid state is represented by a mixture of Fourier series.
		The method is based on splitting the Fokker-Planck equation represented by an integro-partial differential equation into the partial differentiation part for continuous diffusion and the integral part for discrete transition, and integrating the solution of each part.
		The propagated density function is used with a likelihood function in the Bayes' formula to estimate the hybrid state for given sensor measurements.
        The unique feature of the proposed technique is that the probability density describing complete stochastic properties of the hybrid state is constructed without relying on the common Gaussian assumption, in contrast to other methods that compute certain expected values such as mean or variance.
		We apply the proposed method to the estimation of two examples: the bouncing ball model and the Dubins vehicle model. 
		We show that the proposed technique yields propagated densities consistent with Monte Carlo simulations, more accurate estimates over the Gaussian based approach and some computational benefits over a particle filter.
	\end{abstract}
\end{frontmatter}


\section{Introduction}
General Stochastic Hybrid System (GSHS) is a stochastic dynamical systems that incorporates both continuous and discrete dynamics \cite{bujorianu2006toward}.
A GSHS describes the time evolution of a hybrid state, consisting of a discrete state which takes values in a set of countably many discrete modes; and a continuous state which takes values in Euclidean space, or more generally a smooth manifold.
The dynamics of the continuous state is governed by a set of stochastic differential equations (SDEs) associated with the discrete modes; and the discrete jump is usually triggered by the continuous state entering a guard set, or by a Poisson rate function dependent on the hybrid state.
During the discrete jump, the hybrid state is reset according to a stochastic kernel function.

GSHS is one of the most general formalism to inject uncertainties into a hybrid system.
For example, piece-wise deterministic Markov process \cite{davis1984piecewise} requires the evolution of continuous state to be deterministic (i.e., governed by ordinary differential equations), and therefore it can be considered as a special case of GSHS.
Switching diffusion process \cite{ghosh1993optimal} assumes that the discrete jump does not reset the continuous state, so its kernel function can be represented as a delta function over the continuous space, thus it is another special case of GSHS.
The common Markov jump linear system \cite{costa2006discrete} further requires the triggering of the discrete jump does not depend on the continuous state, and the continuous dynamics is linear;
the piecewise affine model is similar, but requires the discrete transitions can be only triggered by polyhedral guard sets.
With its great versatility, GSHS has been widely used to model complicated dynamical systems with uncertainties, such as air traffic management \cite{blom2009rare,prandini2008application,seah2009stochastic}, chemical process \cite{hespanha2005stochastic}, neuron activities \cite{pakdaman2010fluid}, or communicating networks \cite{hespanha2004stochastic}.

In this paper, we study the uncertainty propagation and Bayesian estimation problem for GSHS, i.e., for a given initial probability density function (PDF) of a hybrid state, we wish to compute its PDF at any time in the future and possibly conditioned by noisy measurements.
Previous studies on the state estimation of stochastic hybrid systems mostly focused on some special cases of GSHS.
The model received most attention is the switching diffusion process, and the typical method to estimate its hybrid state is the multiple model approach, which uses a bank of Kalman filters to keep track of the state trajectories.
Various methods have been proposed to deal with the exponentially increasing number of trajectories.
The interacting multiple model (IMM) technique merges trajectories under the same mode by matching a Gaussian distribution to Gaussian mixtures \cite{blom1988interacting}.
Alternatively, only the several most probable trajectories are tracked in \cite{hofbaur2002mode}.
The moving horizon method \cite{ferrari2002moving} optimizes a Kalman type loss function, which is essentially picking out the most probable trajectory.
The method in \cite{zhang1999optimal} averages the dynamics of different modes, leading to a single filter solution.
The Kalman filter in the multiple model method can be replaced by extended or unscented Kalman filters to deal with the non-linearity in the continuous dynamics \cite{prakash2012state}.
However, these methods do not incorporate the reset of continuous state during a discrete jump.
In addition, system identification algorithms have also been proposed for the piecewise affine model to estimate unknown parameters~\cite{bako2011identification,bemporad2005bounded,ferrari2003clustering,juloski2005bayesian,pillonetto2016new,ohlsson2013identification}.
These should be distinguished from the proposed estimation scheme designed for the unknown state of a hybrid system.

Estimation of GSHS with a continuous state reset has been considered in \cite{benazera2009set}, where a set theoretic method was applied to a different formalism of stochastic hybrid system. Thus, it is not readily applicable to GSHS.
In \cite{liu2014hybrid}, the authors studied the continuous filtering problem for GSHS by solving the Zakai equation.
This is achieved by a two-step process: a GSHS is first approximated by a discrete-time Markov chain, and then the discretized Zakai equation is solved using a finite difference method.
Finally, we should mention that the Sampling-Importance-Resampling (SIR) based particle filter is applicable to the most general GSHS \cite{koutsoukos2002monitoring,blom2004particle,tafazoli2006hybrid}.
Nonetheless, although the particles in a particle filter also encode the information of PDF, it is usually computational intensive to recover the PDF from particles, especially when the number of them is large.

Here, we propose a spectral technique for uncertainty propagation and Bayesian estimation of GSHS, which can be regarded as an alternative to the finite difference method adopted in \cite{liu2014hybrid}.
The objective is to construct the probability density function of hybrid state directly, which contains much more stochastic information of the hybrid state than mere the estimates.
According to the framework of Bayesian estimation, it is composed of uncertainty propagation between measurements, and correction when a  measurement becomes available.

The uncertainty propagation part is addressed by solving the Fokker-Planck (FP) equation constructed for GSHS \cite{bec2010unifying,hespanha2004stochastic}, which is an integro-partial differential equation (IPDE) describing the evolution of the PDF.
It is solved by a splitting method: the IPDE is split into the partial differentiation part for continuous diffusion, and the integration part for stochastic discrete transition.
The former is solved by a spectral technique where the PDF is represented by Fourier series, so that its evolution is described by a linear ordinary differential equation (ODE) system of the Fourier coefficients. 
The latter part of IPDE is tackled by using a quadrature rule to approximate the integration term, so that the evolution of PDF is described by another linear ODE system.
The PDF and its Fourier coefficients are converted back and forth to combine the two parts of solutions.
This part of work has been presented in \cite{wang2019}.
 
Next, for the correction part, the posterior PDF is obtained according to Bayes' theorem for given sensor measurements.
In short, we show that the Fokker-Planck equation for uncertainty propagation and estimation can be transformed into a set of ODEs via the spectral technique. 
The superiority of the proposed method is an arbitrary form of the estimated PDF is constructed directly, without limiting to a certain canonical form such as a Gaussian distribution, or the need to count particles in a grid.
As such, this should be distinguished from other estimation techniques where only selected stochastic properties, such as mean and variance, are computed.

To verify the validity and versatility of the proposed uncertainty propagation and estimation method, we apply it to two examples: the bouncing ball model and the Dubins vehicle model.
Despite its simplicity (only one discrete mode), the bouncing ball model is a typical GSHS that cannot be modeled as a switching diffusion process, because it resets the velocity of the ball during its collision with the ground.
The Dubins vehicle model is a more complex one, with three discrete modes and three dimensional continuous space, and is used to represent a scenario in air traffic control where the state of an airplane avoiding collisions with obstacles is estimated.
We first verify the propagated uncertainty without measurements using the proposed method is consistent with Monte Carlo simulations.
Then with noisy measurements, we compare the proposed technique with a particle filter in both examples, and with an IMM in the Dubins vehicle model.
Results indicate the proposed method is superior both in estimation accuracy and computation efficiency than the IMM method, and has similar accuracy but faster computation speed compared to the particular particle filter implemented.
The MATLAB implementation can be found in \cite{HybridCode}.

The remainder of this paper is organized as follows.
In Section \ref{section:problem}, GSHS and its FP equation are introduced, and the Bayesian estimation problem is formulated.
In Section \ref{section:method}, the numerical method is developed to solve the FP equation, and the resulting PDF is used in Bayes' formula to estimate the hybrid state.
This method is applied to the bouncing ball model and Dubins vehicle model in Section \ref{section:ball} and \ref{section:car}, followed by discussions and conclusions in Section \ref{sec:discussion} and \ref{section:conclusion}.


\section{Problem Formulation} \label{section:problem}
In this section, we introduce the definition of GSHS and the corresponding Fokker-Planck equation that describes the evolution of its PDF.  Later, we formulate the Bayesian estimation problem considered in this paper.

\subsection{General Stochastic Hybrid System (GSHS)}
The GSHS is defined as a collection $H=\left\{X,a,b,Init,\lambda,K\right\}$ described as follows:
\begin{itemize}
	\item $X=\left(\mathbb{R}^{N_r},S\right)$ is the hybrid state space, where $\mathbb{R}^{N_r}$ is $N_r$-dimensional Euclidean space, and $S$ is the set of $N_s$ discrete modes.
	The hybrid state is denoted by $x=(r,s)\in X$, where $r\in\mathbb{R}^{N_r}$ is the continuous state, and $s\in S$ is the discrete mode.
	\item $Init: \mathcal{B}(X)\to[0,1]$ is the initial distribution of the hybrid state, where $\mathcal{B}(X)$ is all Borel sets in $X$.
	\item The continuous state evolves according to the following stochastic differential equation:
		\begin{equation} \label{eqn:ContSDE}
			\text{d}r = a(t,r,s)\text{d}t + b(t,r,s)\text{d}W_t,
		\end{equation}
	where $a:(\mathbb{R}^+,X)\to\mathbb{R}^{N_r}$ is the drift vector, $b:(\mathbb{R}^+,X)\to\mathbb{R}^{N_r\times N_w}$ is the coefficient matrix of diffusion, and $W_t$ is a standard $N_w$-dimensional Wiener process.
	\item The discrete transition is triggered by a Poisson rate function, $\lambda:X\to\mathbb{R}^+$, i.e. the probability that a transition happens in a time interval $\Delta t$ is $\lambda(r,s)\Delta t + o(\Delta t)$.
	\item During each discrete transition, the hybrid state is reset according to the stochastic kernel $K:X\times\mathcal{B}(X)\to[0,1]$. Specifically, $K(x^-,X^+)$ gives the probability that the hybrid state is reset into $X^+\subseteq X$ from $x^-$.
\end{itemize}

A notable restriction of the GSHS considered here compared to \cite{bujorianu2006toward} is that we do not consider the ``forced'' discrete jump which is triggered by the continuous state entering a guard set.
However, the forced jump can be approximated by a ``spontaneous'' jump triggered by a rate function as in \cite{abate2008approximation,hespanha2004stochastic}.
Specifically, if the rate function is set to zero outside the guard set, and it is set to be large inside the guard set, then the jump triggered by the rate function closely mimics the jump triggered by the guard set.
Later, we demonstrate this approximation in the bouncing ball example, where the discrete jump happens when the ball touches the ground with a negative velocity.

Let $(\Omega,\mathcal{F},\mathbb{P})$ be the underlying probability space, where $\Omega$ is the sample space, $\mathcal{F}$ is a $\sigma$-algebra over $\Omega$, and $\mathbb{P}$ denotes the probability measure.
We further assume the distribution of the hybrid state has a PDF $p(t,r,s)$ such that for all measurable $A\subseteq X$, $\mathbb{P}\{(r,s)\in A\} = \sum_{s\in S}\int_{\{r:\;(r,s)\in A\}}p(t,r,s)dr$.
Also, the kernel function can be represented as a set of density functions: $\kappa(x^-,x^+)$, such that $K(x^-,X^+)=\int_{X^+}\kappa(x^-,x^+)\text{d}x^+$.

For a given $\omega\in\Omega$, let $\{u_k(\omega)\}$ be a sequence of independent, uniform random variables on $[0,1]$.
An execution (a right-continuous sample path) of a GSHS can be generated as follows:
\begin{enumerate}
	\item Initialize $r(\omega,t_0) = r_0$, $s(\omega,t_0) = s_0$ according to the initial distribution.
	\item \label{second step} Let $t_1(\omega) = \text{sup}\left\{t:e^{-\int_{t_0}^{t}\lambda(r(\omega,t^{\prime}),s(\omega,t^{\prime}))\text{d}t^{\prime}}>u_1(\omega)\right\}$ be the first jump time.
	\item During $t\in[t_0,t_1(\omega))$, $s(\omega,t)=s_0$ and $r(\omega,t)$ is a sample path of SDE (\ref{eqn:ContSDE}) with $s=s_0$.
	\item At time $t_1$, the state is reset to $r(\omega,t_1^+)$ and $s(\omega,t_1^+)$ as a sample from the kernel $\kappa(r(\omega,t_1^-),s_0,r^+,s^+)$.
	\item if $t_1<\infty$, repeat from \ref{second step}) with $t_0$, $s_0$, $t_1$, $u_1$ replaced by $t_k(\omega)$, $s(\omega,t_k^+)$, $t_{k+1}(\omega)$, $u_{k+1}$ for $k=1,2,\ldots$
\end{enumerate}

\subsection{Uncertainty Propagation and Estimation}

For a given GSHS and its initial PDF $p(0,r,s)$, we wish to determine the PDF of the hybrid state at any time in the future, conditioned by the measurements available.
The resulting PDF can be used to compute any stochastic property of the estimated hybrid state. 
According to the Bayesian estimation formulation, it is composed of the uncertainty propagation and the correction.

First, the propagation of the PDF can be calculated by solving the following Fokker-Planck equation generalized for GSHS \cite{hespanha2004stochastic}, 
\begin{align} \label{eqn:FP} 
	\frac{\partial p(t,r,s)}{\partial t} = &-\sum_{i=1}^{N_r}\frac{\partial}{\partial r_i}\left(a_i(t,r,s)p(t,r,s)\right) \nonumber\\
	&+ \sum_{i,j=1}^{N_r}\frac{\partial^2}{\partial r_i\partial r_j}\left(D_{ij}(t,r,s)p(t,r,s)\right) \nonumber \\
	&+ \sum_{s^{\prime}\in S}\int_{\mathbb{R}^{N_r}}\kappa(r^{\prime},s^{\prime},r,s)\lambda(r^{\prime},s^{\prime})p(t,r^{\prime},s^{\prime})\text{d}r^{\prime} \nonumber \\
	&- \lambda(r,s)p(t,r,s),
\end{align}
where $r_i$, $a_i$ are the $i$-th component of $r$ and $a$ respectively, $D=\frac{1}{2}bb^T$, and $D_{ij}$ is the $i,j$-th component of $D$.
Note that (\ref{eqn:FP}) represents a set of $N_s$ IPDEs as $s$ varies in $S$.
For a specific $s$, (\ref{eqn:FP}) can be decomposed into the contribution due to the continuous SDE and the discrete jump, corresponding to the first two terms, and the last two terms on the right-hand side (RHS) of (\ref{eqn:FP}), respectively.
To clearly distinguish these two parts, it is rewritten as:
\begin{equation} \label{eqn:FPSplitted}
	\frac{\partial p(t,r,s)}{\partial t} = \mathcal{L}_c^*p(t,r,s) + \mathcal{L}_d^*p(t,r,s),
\end{equation}
where $\mathcal{L}_c^*$ and $\mathcal{L}_d^*$ denote the adjoint of the infinitesimal generators of the continuous SDE and the discrete jump respectively.
Because of the presence of both partial differentiation and integration in (\ref{eqn:FP}), solving it is fundamentally challenging.
In the next section, we present a numerical method to tackle this problem based on the splitting and spectral methods.

Second, for the correction part, we assume the hybrid state is observed by a measurement function
\begin{equation}
	z_k = h_k(r_k,s_k,v_k),
\end{equation}
at a sequence of time $\{t_k\}$, where $v_k\in V$ is a sequence of independent measurement noises in arbitrary noise space with density function $p_{v,k}(v)$.
There is no general way to calculate the likelihood function from $h$ and $p_v$, but when the noise is addictive, i.e., $z_k=h_k(r_k,s_k)+v_k$, the likelihood function has a simple form as $p_{v,k}(z_k-h_k(r_k,s_k))$.
Here, we assume the likelihood function is given as $p_k(z_k\lvert x_k)$, the goal of this paper is to find the posterior PDF $p\left(x_k\lvert z_{1:k}\right)$.
The posterior PDF can be easily distilled into any estimation results, such as the mean, variance, mode of the continuous state, and the most probable discrete state, etc.


\section{Spectral Bayesian Estimation} \label{section:method}

The proposed method solves the FP equation \eqref{eqn:FPSplitted} by splitting it into two equations, solving them independently and combining them through the splitting method \cite{florchinger1991time}:
\begin{subequations}
	\begin{align}
	\frac{\partial p^c(t,r,s)}{\partial t} = \mathcal{L}_c^*p^c(t,r,s) \label{eqn:FPCont} \\
	\frac{\partial p^d(t,r,s)}{\partial t} = \mathcal{L}_d^*p^d(t,r,s) \label{eqn:FPDisc}
	\end{align}
\end{subequations}
The first equation \eqref{eqn:FPCont} caused by continuous diffusion is a PDE, and we solve it using the spectral method; the second equation \eqref{eqn:FPDisc} caused by discrete jump is an integral equation, and we approximate it with an ODE system using a quadrature rule.
The resulting uncertainty propagation technique is used with a correction step to conduct Bayesian estimation of the GSHS.

\subsection{Spectral Method for Continuous Evolution} \label{section:Spectral}

It is well known that the Fourier series $\left\{e^{2\pi inx/L}\right\}_{n=-\infty}^{\infty}$ form an orthonormal basis of $\mathcal{L}^2([-\frac{L}{2},\frac{L}{2}])$, i.e. the space of square integrable functions on $[-\frac{L}{2},\frac{L}{2}]$, for any $L>0$ \cite{stein2009real}.
In addition, the PDF of a real valued random variable $p\in\mathcal{L}^1(\mathbb{R})$ is the limit of a series of compact supported functions $\{p_n\}$.
Let $p_n$ be a ``close enough" approximation of $p$ in $\mathcal{L}^1$ sense.
Then, $p_n$ can be regarded as a function in $\mathcal{L}^2([-\frac{L}{2},\frac{L}{2}])$ for some $L$, as it has a compact support and is integrable.
As a result, $p_n$ can be further represented by its Fourier series.
The above argument can be extended to the PDF of a GSHS, as $p(t,r,s)\in\mathcal{L}^1(\mathbb{R}^{N_r})$ for any $s\in S$ and $t\in\mathbb{R}^+$.
Therefore, in this section we assume $p(t,r,s)$ can be approximated to an arbitrary precision with its truncated Fourier series on $r$, for any time $t$ and discrete state $s$.

The benefit of using Fourier series is that the PDF is not restricted to a specific form, such as the Gaussian distribution that is commonly used in assumed density filters.
Furthermore, the concept of Fourier series can be extended to a compact Lie group, such as the rotational group SO(3), by Peter-Weyl Theorem \cite{PetWeyMA27} which states that the matrix components of the irreducible unitary representations form a complete orthonormal basis for the space of square integrable functions defined on that group.
This means the numerical methods developed in this section can be extended to a GSHS with a compact Lie group as the continuous state space, which can be very useful if the attitude of a rigid body is considered.
We will study a simpler example in the Dubins vehicle model where the heading angle lies in the circular space homeomorphic to the unit circle $\mathbb{S}^1$.

In this subsection, we use the spectral method to solve \eqref{eqn:FPCont} caused by the continuous SDE.
Equation \eqref{eqn:FPDisc} caused by discrete jump is considered in Section \ref{section:quadrature}.
The key idea is to use discrete Fourier transform (DFT) to convert this PDE into a linear ODE system which can be solved efficiently using matrix exponential.
To derive this transformation in a concise manner, we first consider a simpler case where the dimension of the continuous space is one, and at the end of this sub-section we show how it can be generalized to multi-dimensional cases.

Let $\hat{f}_n[g]$ denote the $n$-th Fourier coefficient of a function $g\in\mathcal{L}^2([-\frac{L}{2},\frac{L}{2}])$, where $n\in\left\{-\frac{N}{2},-\frac{N}{2}+1,\ldots,\frac{N}{2}-1\right\}$ for some even number $N$.
The Fourier coefficient $\hat{f}_n[g]$ can be calculated using DFT:
\begin{equation} \label{eqn:DFT}
	\hat{f}_n[g] = \frac{1}{N}\sum_{j=-\frac{N}{2}}^{\frac{N}{2}-1}g(r_j)e^{-\frac{2\pi inr_j}{L}},
\end{equation}
where $r_j=\frac{jL}{N}$ is a pre-defined grid on $[-\frac{L}{2},\frac{L}{2})$ with $N$ points.
The approximation of function $g$ can be recovered using inverse DFT (IDFT):
\begin{equation} \label{eqn:IDFT}
	\tilde{g}(r) = \sum_{n=-\frac{N}{2}}^{\frac{N}{2}-1}\hat{f}_n[g]e^{\frac{2\pi inr}{L}}.
\end{equation}
Equation (\ref{eqn:IDFT}) is exact at the grid points $r_j$, and furthermore it is a good interpolation of function $g$ between two grid points \cite{shen2011spectral}.
In the following development, we do not distinguish $g$ and $\tilde{g}$.

Another desirable feature of the Fourier transform is that there are several operational properties, where the Fourier coefficients of differentiation and product of functions can be represented by their own Fourier coefficients.
The formulae for DFT-based differentiation in frequency domain can be found in \cite{joh2011DFT}:
\begin{align} \label{eqn:DFTDiff}
	\begin{split}
		\hat{f}_n\left[\frac{\text{d}g}{\text{d}r}\right] &= \frac{2\pi inc_n}{L}\hat{f}_n[g], \\
		\hat{f}_n\left[\frac{\text{d}^2g}{\text{d}r^2}\right] &= -\frac{4\pi^2n^2}{L^2}\hat{f}_n[g].
	\end{split}
\end{align}
where $c_n=0$ for $n=-\frac{N}{2}$ and $c_n=1$ otherwise.
We can also use the discrete convolution to calculate the Fourier coefficients of the product of two functions:
\begin{equation} \label{eqn:DFTProd}
	\hat{f}_n[gh] = \sum_{k=-\frac{N}{2}}^{\frac{N}{2}-1}\hat{f}_{n-k}[g]\hat{f}_k[h].
\end{equation}
Using (\ref{eqn:DFTDiff}) and (\ref{eqn:DFTProd}), the Fourier coefficients of the RHS of (\ref{eqn:FPCont}) can be expressed as a linear combination of the Fourier coefficients of $p$ on $r$, given any time $t$ and discrete state $s$, which is summarized in the following theorem.

\begin{thm}
	The Fourier coefficients of the PDF $p^c(t,r,s)$ in (\ref{eqn:FPCont}) evolve approximately according to $N_s$ linear ODE systems corresponding to all discrete modes $s\in S$:
	\begin{equation} \label{eqn:FPContODE}
		\frac{\mathrm{d}\hat{f}[p^c](t,s)}{\mathrm{d}t} = A(t,s)\hat{f}[p^c](t,s),
	\end{equation}
	where $\hat{f}[p^c](t,s) = \left[\hat{f}_{-\frac{N}{2}}[p^c](t,s), \ldots, \hat{f}_{\frac{N}{2}-1}[p^c](t,s)\right]^T$, and $\hat{f}_n[\cdot](t,s)$ means the Fourier coefficient is dependent on time $t$ and discrete mode $s$.
	Also, the matrix $A(t,s)\in\mathbb{C}^{N\times N}$ is defined as:
	\begin{equation} \label{eqn:ContMatrix}
		A(t,s) = \left[\begin{matrix}
			d_{-\frac{N}{2},-\frac{N}{2}} & \ldots & d_{-\frac{N}{2},\frac{N}{2}-1} \\
			\vdots & \ddots & \vdots \\
			d_{\frac{N}{2}-1,-\frac{N}{2}} & \ldots & d_{\frac{N}{2}-1,\frac{N}{2}-1}
		\end{matrix}\right],
	\end{equation}
	where
	\begin{equation} \label{eqn:contCoef}
		d_{j,k}=-\frac{2\pi ijc_j}{L}\hat{f}_{j-k}[a](t,s) 
		-\frac{4\pi^2j^2}{L^2}\hat{f}_{j-k}[D](t,s).
	\end{equation}
\end{thm}
\textit{Proof}: Denote the RHS of (\ref{eqn:FPCont}) by $F$, represent both sides by their Fourier series, we get
\begin{equation}
	\frac{\partial}{\partial t}\sum_{n=-\frac{N}{2}}^{\frac{N}{2}-1}\hat{f}_n[p^c](t,s)e^{\frac{2\pi inr}{L}} = \sum_{n=-\frac{N}{2}}^{\frac{N}{2}-1}\hat{f}_n[F](t,s)e^{\frac{2\pi inr}{L}}.
\end{equation}
Because the Fourier series are orthonormal basis, each of the coefficients of $e^{2\pi inr/L}$ can be equated, so $\frac{\text{d}}{\text{d}t}\hat{f}_n[p^c](t,s) = \hat{f}_n[F](t,s)$.
By noting the two terms in $F$ are the products of $p$ with known functions $a$ and $D$ (one dimensional here), followed by differentiation operations, the product formula (\ref{eqn:DFTProd}) and differentiation formulae (\ref{eqn:DFTDiff}) can be used to express $\hat{f}_n[F](t,s)$ in terms of $\hat{f}_n[p](t,s)$, then we get the following ODE system:
\begin{align} \label{eqn:FPContDFT}
	\frac{\text{d}}{\text{d}t}\hat{f}_n[p^c](t,s) = &-\frac{2\pi inc_n}{L}\sum_{k=-\frac{N}{2}}^{\frac{N}{2}-1}\hat{f}_{n-k}[a](t,s)\hat{f}_k[p^c](t,s) \nonumber \\
	&-\frac{4\pi^2n^2}{L^2}\sum_{k=-\frac{N}{2}}^{\frac{N}{2}-1}\hat{f}_{n-k}[D](t,s)\hat{f}_k[p^c](t,s),
\end{align}
from which $A(t,s)$ can be derived. \hspace*{\fill} $\square$

Note that the Fourier coefficient is $N$-periodic, i.e. $\hat{f}_{n+kN}[\cdot]=\hat{f}_{n}[\cdot]$ for any integer $k$.
So in (\ref{eqn:contCoef}), if $j-k<-\frac{N}{2}$ or $j-k>\frac{N}{2}-1$, the index can be wrapped into $[-\frac{N}{2},\frac{N}{2}-1]$ based on this peroidicity.
In addition, if the SDE of the continuous state (\ref{eqn:ContSDE}) is time-invariant, i.e. $a(r,s)$ and $b(r,s)$ do not depend on time $t$, then the coefficient matrices $A(s)$ also do not depend on $t$, and the ODE systems become time-invariant and can be analytically solved by matrix exponential:
\begin{equation} \label{eqn:ContSol}
	\hat{f}[p^c](t,s) = e^{A(s)t}\hat{f}[p^c](0,s).
\end{equation}
The numerical algorithms to calculate matrix exponential have been summarized and compared in \cite{moler2003nineteen}.
Though this may be computational intensive for large matrices in higher dimensional space, they only need to be calculated once in implementation, and the only repeated calculation is matrix-vector multiplication.

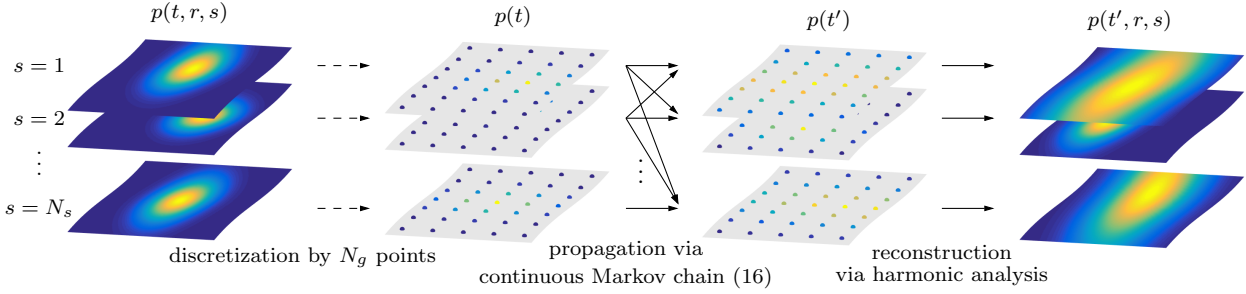
\begin{figure*}[h]
	\setlength{\unitlength}{2.05em}\footnotesize
	\centerline{
		\includegraphics[width=0.9\textwidth]{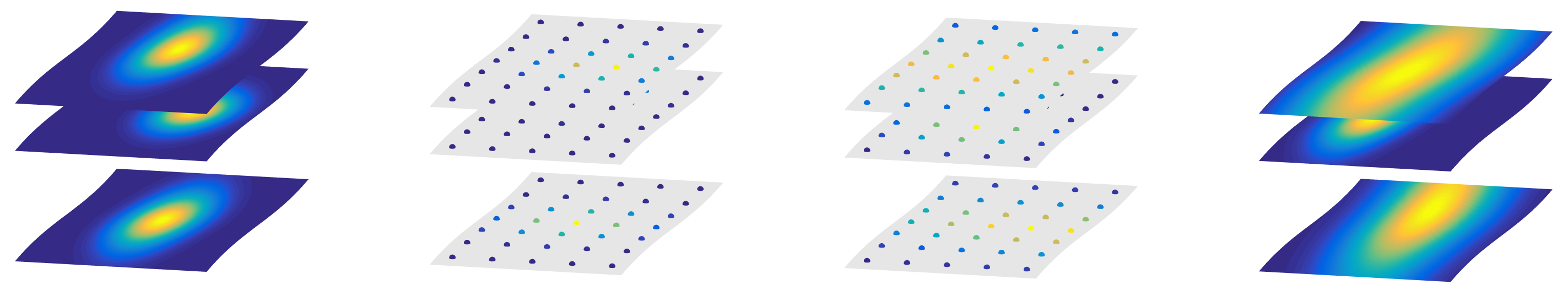}
	}
	\centerline{
		\begin{tikzpicture}[overlay,shift={(-8.2,0)}]
		\draw[dashed,arrows={-Triangle[angle=30:4pt]}] (3.7,2.9) -- ++(0:0.7);
		\draw[dashed,arrows={-Triangle[angle=30:4pt]}] (3.7,2.2) -- ++(0:0.7);
		\draw[dashed,arrows={-Triangle[angle=30:4pt]}] (3.7,1.0) -- ++(0:0.7);
		\node at (3.5,0.35) {\scriptsize discretization by $N_g$ points};
		\draw[arrows={-Triangle[angle=30:4pt]}] (7.8,2.9) -- ++(0:0.7);
		\draw[arrows={-Triangle[angle=30:4pt]}] (7.8,2.9) -- ++(0.7,-0.65);
		\draw[arrows={-Triangle[angle=30:4pt]}] (7.8,2.9) -- ++(0.7,-1.85);
		\draw[arrows={-Triangle[angle=30:4pt]}] (7.8,2.2) -- ++(0:0.7);
		\draw[arrows={-Triangle[angle=30:4pt]}] (7.8,2.2) -- ++(0.7,0.65);
		\draw[arrows={-Triangle[angle=30:4pt]}] (7.8,2.2) -- ++(0.7,-1.15);
		\draw[arrows={-Triangle[angle=30:4pt]}] (7.8,1.0) -- ++(0:0.7);
		\node at (8.0,1.6) {$\vdots$};
		\node at (7.8,0.25) {\scriptsize \shortstack[c]{propagation via\\ continuous Markov chain \eqref{eqn:FPDiscODE}}};
		\draw[arrows={-Triangle[angle=30:4pt]}] (12.0,2.9) -- ++(0:0.7);
		\draw[arrows={-Triangle[angle=30:4pt]}] (12.0,2.2) -- ++(0:0.7);
		\draw[arrows={-Triangle[angle=30:4pt]}] (12.0,1.0) -- ++(0:0.7);
		\node at (12.0,0.25) {\scriptsize \shortstack[c]{reconstruction\\ via harmonic analysis}};
		\node at (0.0,2.9) {\scriptsize $s=1$};
		\node at (0.0,2.2) {\scriptsize $s=2$};
		\node at (0.0,1.7) {\scriptsize $\vdots$};
		\node at (0.0,1.0) {\scriptsize $s=N_s$};
		\node at (2.0,3.6) {\scriptsize $p(t,r,s)$};
		\node at (6.3,3.55) {\scriptsize $p(t)$};
		\node at (10.5,3.5) {\scriptsize $p(t')$};
		\node at (14.5,3.5) {\scriptsize $p(t',r,s)$};
		\end{tikzpicture}
	}
	\vspace*{0.2cm}
	\caption{Uncertainty propagation through discrete jumps: the sample values of the probability density are obtained at $\{r_i\}_{i=1}^{N_g}$ for each $s\in S$.
	Then, the Fokker--Planck equation describing the evolution of a probability density over discrete transitions is approximated by a continuous-time Markov chain and the sample values of the probability density are propagated via \eqref{eqn:DiscSol} explicitly, from which the propagated density function is reconstructed using  harmonic analysis.}\label{fig:UPD}
\end{figure*}

Finally we show how the spectral transform approach can be generalized to multi-dimensional continuous space.
Suppose $r\in\mathbb{R}^{N_r}$, the Fourier series can be replaced by $\text{exp}\left\{2\pi i\left(\frac{n_1r_1}{L_1}+\frac{n_2r_1}{L_2}+\ldots+\frac{n_{N_r}r_{N_r}}{L_{N_r}}\right)\right\}$.
The summation in the DFT (\ref{eqn:DFT}) and IDFT (\ref{eqn:IDFT}) is changed to multi-summation along each axis, with the grid defined as $r_{j_1j_2\ldots j_{N_r}}=\left(\frac{j_1L_1}{N_1},\frac{j_2L_2}{N_2},\ldots,\frac{j_{N_r}L_{N_r}}{N_{N_r}}\right)$, where $j_k\in\left\{-\frac{N_k}{2},\ldots,\frac{N_k}{2}-1\right\}$ for $1\leq k\leq N_r$.
The equations (\ref{eqn:DFTDiff}) do not change if the differentiation is replaced by partial differentiation.
The product formula (\ref{eqn:DFTProd}) remains true if the RHS is replaced by multi-dimensional discrete convolution.
From these, the ODE system (\ref{eqn:FPContDFT}) can be similarly derived for multi-dimensional continuous space, which is given in (\ref{eqn:FPContDFTMulti}).

\subsection{Propagation Over Discrete Transition } \label{section:quadrature}

Next, we consider the evolution of PDF due to the stochastic discrete transition, as given by \eqref{eqn:FPDisc}.
For notation simplicity, let the grid points over the continuous state space defined in the last paragraph of Section \ref{section:Spectral} be re-parameterized by $\left\{r_j\right\}$, where $j\in\left\{1,2,\ldots,N_g\right\}$, and $N_g=\prod_{k=1}^{N_r}N_k$ is the total number of grid points.
The key idea of the proposed approach is that the evolution of the PDF on the grid over \eqref{eqn:FPDisc} is described by a linear ODE system as follows.
\begin{thm}
	The values of the PDF $p^d(t,r_j,s)$ on the grid points in (\ref{eqn:FPDisc}) evolve approximately according to a linear ODE system:
	\begin{equation} \label{eqn:FPDiscODE}
		\frac{\mathrm{d}p^d(t)}{\mathrm{d}t} = Bp^d(t),
	\end{equation}
	where $p^d(t)$ is all $\left[p^d(t,r_1,s),\ldots,p^d(t,r_{N_g},s)\right]^T$ concatenated into a column vector from $s=s_1$ to $s=s_{N_s}$.
	The coefficient matrix $B$ is given in (\ref{eqn:DiscMatrix}).
\end{thm}
\textit{Proof}: Replace the integral term in (\ref{eqn:FPDisc}) by a finite summation, it can be approximated by:
\begin{align} \label{eqn:FPDiscNum} 
	&\frac{\partial p^d(t,r_i,s)}{\partial t} = -\lambda(r_i,s)p^d(t,r_i,s) + \nonumber\\
	& \quad \sum_{s^{\prime}\in S}\sum_{j=1}^{N_g}\kappa(r_j^{\prime},s^{\prime},r_i,s)\lambda(r_j^{\prime},s^{\prime})p^d(t,r_j^{\prime},s^{\prime})\Delta(r^{\prime}_j),
\end{align}
where $\Delta(r^{\prime}_j)$ is the weight for a quadrature rule.
Because $\kappa$ and $\lambda$ are fixed functions and do not depend on $t$, (\ref{eqn:FPDiscNum}) is a time invariant linear ODE system of $p^d(t,r_j,s)$, for $j=1,\ldots,N_g$ and $s=1,\ldots,N_s$.
And the coefficient matrix $B$ can be derived from (\ref{eqn:FPDiscNum}). \hspace*{\fill} $\square$

The simplest form of weight $\Delta(r^{\prime}_j)$ in (\ref{eqn:FPDiscNum}) is $\prod_{k}^{N_r}L_k/N_k$, i.e. the volume of one grid block.
Similar to (\ref{eqn:ContSol}), the solution to (\ref{eqn:FPDiscNum}) can also be solved analytically using matrix exponential:
\begin{equation} \label{eqn:DiscSol}
	p^d(t) = e^{Bt}p^d(0),
\end{equation}
Once the probability density at each grid point is propagated via \eqref{eqn:DiscSol}, the Fourier series can be calculated from DFT (\ref{eqn:DFT}), and the PDF between grid points can be interpolated by IDFT (\ref{eqn:IDFT}). 
These are illustrated in Fig. \ref{fig:UPD}.

In general, the matrix $B$ is a $\prod_{k=1}^{N_r}N_kN_s$-by-$\prod_{k=1}^{N_r}N_kN_s$ real valued matrix, but for a switching diffusion process, (\ref{eqn:FPDiscNum}) can be significantly simplified as follows.
Because in a switching diffusion process, the continuous state is not reset, the kernel function can be written with a Dirac delta function over the continuous state:
\begin{equation} \label{eqn:KernelSDP}
	\kappa(r^-,s^-,r^+,s^+) = \delta(r^+-r^-)\tilde{\kappa}(r^-,s^-,s^+).
\end{equation}
Substitute (\ref{eqn:KernelSDP}) into (\ref{eqn:FPDisc}), the integral over delta function can be simplified:
\begin{align} \label{eqn:FPDiscNumSDP}
	\frac{\partial p^d(t,r,s)}{\partial t} = &\sum_{s^-\in S}\tilde{\kappa}(r,s^-,s)\lambda(r,s^-)p^d(t,r,s^-) \nonumber \\
	&- \lambda(r,s)p^d(t,r,s).
\end{align}
Equation (\ref{eqn:FPDiscNumSDP}) is a linear ODE system of $p^d(t,r,s)$ for $s=1,\ldots,N_s$ at any $r\in\mathbb{R}^{N_r}$, so the off-diagonal terms in $B$ between different grid points $r_i$ and $r_j$ are eliminated.
Taking values at those grid points, (\ref{eqn:FPDiscNumSDP}) becomes a set of $N_g$ ODE systems, each of which is of dimension $N_s$.
Intuitively, since the discrete jump does not reset the continuous state, the marginal distribution of the continuous state does not change, and what happens during the jump is the point masses of the discrete modes switch between each other at each $r$.

\subsection{Splitting Method} \label{section:splitting}

\begin{figure}
	\setlength{\unitlength}{2.2em}\footnotesize
	\centerline{
		\begin{picture}(6,4.5)(0.0,-3.1)
		\put(0,0){\framebox(2,1.0)[c]{\shortstack[c]{$p_k^c=p_k$}}}
		\put(4.2,0){\framebox(2.3,1.0)[c]{\shortstack[c]{$p_{k+1}^c$}}}
		\put(0,-2){\framebox(2,1.0)[c]{\shortstack[c]{$p_k^d=p_{k+1}^c$}}}
		\put(4.2,-2){\framebox(2.3,1.0)[c]{\shortstack[c]{$p_{k+1}=p_{k+1}^d$}}}
		\put(2,0.5){\vector(1,0){2.2}}
		\put(2,-1.5){\vector(1,0){2.2}}
		\put(4.1,0){\vector(-2,-1){2}}
		\put(0.5,-2.8){$t=t_{k}$}
		\put(4.65,-2.8){$t=t_{k+1}$}
		\put(2.2,0.4){\parbox{0.6in}{\centering cont. dyn. (\ref{eqn:FPCont})}}
		\put(3.4,-0.65){re-initialization}
		\put(2.2,-1.65){\parbox{0.6in}{\centering disc. dyn. (\ref{eqn:FPDisc})}}
		\end{picture}
	}
	\caption{\textit{Splitting methods:} the initial uncertainty $p_k$ is propagated for the continuous process to obtain $p_{k+1}^c$, which is used as the initial condition $p_{k}^d$ of the discrete dynamics to find $p_{k+1}$.}\label{fig:SPM}
\end{figure}
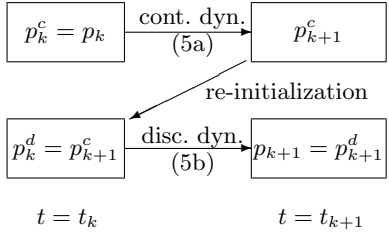

The splitting technique is used to combine the solution of (\ref{eqn:FPCont}) solved in Section \ref{section:Spectral} and the solution of (\ref{eqn:FPDisc}) solved in Section \ref{section:quadrature} as a solution of (\ref{eqn:FPSplitted}).
Let time $t$ be discretized by $\{t_k\}_{k=1}^{\infty}$.
Given a PDF $p(t_k,r,s)$ at $t=t_k$, we initialize $p^c(t_k,r,s)$ as $p(t_k,r,s)$, and solve (\ref{eqn:FPCont}) from $t_k$ to $t_{k+1}$ to get $p^c(t_{k+1},r,s)$.
This is chosen as the initial condition for (\ref{eqn:FPDisc}) at time $t_k$ (re-initialization), i.e. $p^d(t_k,r,s)=p^c(t_{k+1},r,s)$.
Then we solve (\ref{eqn:FPDisc}) again from $t_k$ to $t_{k+1}$ and get $p^d(t_{k+1},r,s)=p(t_{k+1},r,s)$.
See Fig. \ref{fig:SPM}.
This is repeated at each step, and the combined solution converges to the solution to (\ref{eqn:FPSplitted}) if the time-step is sufficiently small \cite{florchinger1991time}.
Note that the values of PDF at grid points and its Fourier coefficients contain equivalent information, and can be converted to each other using DFT (\ref{eqn:DFT}) and IDFT (\ref{eqn:IDFT}).
These two formulae are repeatedly used in the switching between propagating the continuous and discrete dynamics.

It should be noted that the entire equation (\ref{eqn:FP}) can be converted to an ODE system by DFT without resorting to the splitting method \cite{wang2019}.
However, the resulting coefficient matrix of the ODE system is very complicated and hard to compute due to the complexity in the third term on the RHS.

\subsection{Bayesian Estimation}
Next we consider the correction step of Bayesian estimation when a new measurement is available.
With the propagated PDF of the GSHS and the likelihood function, it is straightforward to use the Bayes' formula to calculate the posterior PDF:
\begin{equation} \label{eqn:Bayes}
	p(t_k,r,s \lvert z_{1:k}) \propto p\left(z_k \lvert r,s\right)p(t_k,r,s \lvert z_{1:k-1}),
\end{equation}
where $p(t_k,r,s \lvert z_{1:k-1})$ is the prior PDF propagated from $p(t_{k-1},r,s \lvert z_{1:k-1})$ using the method introduced in this section.
From the posterior distribution, the hybrid state can be estimated as a point estimation problem by minimizing a loss function.
For example, minimizing least squared error leads to the expected value as the estimate, which is typically used for continuous uni-modal distributions.
For multi-modal distributions and the discrete state, using the state that maximizes the posterior PDF is the easiest yet reasonable choice.
The pseudo-code of proposed estimation scheme is listed in Table \ref{tab:EstProc}.

\begin{table}
	\caption{Estimation for general stochastic hybrid system \label{tab:EstProc}}
	\begin{algorithmic}[1]
		\algrule[0.8pt]
		\Procedure{Estimate the hybrid state for GSHS}{}
		\algrule
		\State Initial PDF $p_0$
		\Repeat
		\State $p_{k+1}=${\fontshape{sc}\selectfont Propagation}$(p_k,a_k,b_k)$
		\State $k=k+1$
		\Until{${z}_{k+1}$ is available}
		\State $p_{k+1}=${\fontshape{sc}\selectfont Correction}$(p_{k+1},{z}_{k+1})$
		\State Calculate estimates or other stochastic properties of the state using $p_{k+1}$
		\State \textbf{go to} Step 3
		\EndProcedure
		\algrule
		\Procedure{$p_{k+1}$=Propagation}{$p_k,a_k,b_k$}
		\State Let $p^c_k=p_k$, do DFT and get $\hat{f}[p^c_k]$
		\State Propagate $\hat{f}[p^c_k]$ over \eqref{eqn:ContSol} for time-invariant systems, or solve \eqref{eqn:FPContODE} for time-variant systems, and get $\hat{f}[p^c_{k+1}]$
		\State Do IDFT to $\hat{f}[p^c_{k+1}]$ and get $p^c_{k+1}$
		\State Let $p^d_k=p^c_{k+1}$, propagate $p^d_k$ over \eqref{eqn:DiscSol} and get $p^d_{k+1}$, let $p_{k+1}=p^d_{k+1}$
		\EndProcedure
		\algrule
		\Procedure{$p^+$=Correction}{$p^-,z$}
		\State Use \eqref{eqn:Bayes}, where $z_k=z$, $p(t_k,r,s \lvert z_{1:k-1})=p^-$, and $p(t_k,r,s \lvert z_{1:k})=p^+$
		\EndProcedure
		\algrule[0.8pt]
	\end{algorithmic}
\end{table}


\section{Application to Bouncing Ball} \label{section:ball}

In this section, we apply the proposed uncertainty propagation and estimation methods presented in last section to the classic bouncing ball model.
As noted in the introduction, the bouncing ball model is a representative GSHS as it cannot be modeled as a switching diffusion process that has been repeatedly studied in estimation of stochastic hybrid systems. 

\subsection{Bouncing Ball Model}

In the formulation of GSHS, the bouncing ball model has only one discrete mode, and its continuous state is $r=(y,\dot{y})\in\mathbb{R}^+\times\mathbb{R}$, which is the position (height measured from the ground surface) and velocity of the ball.
The continuous SDE is given by
\begin{equation} \label{eqn:SDEBall}
	\text{d}\left[\begin{matrix}y\\\dot{y}\end{matrix}\right] = 
    \left[\begin{matrix} \dot{y}\\ -g-\nu\dot{y} |\dot{y}| \end{matrix}\right] \text{d}t
	+\left[\begin{matrix}0\\\sigma_v\dot{y}^2\end{matrix}\right]\text{d}W_t,
\end{equation}
where the constant $g\in\mathbb{R}$ is the gravitational acceleration, and $\nu\in\mathbb{R}$ is an atmospheric drag coefficient which is disturbed by a Gaussian white noise with strength $\sigma_v\in\mathbb{R}$.
In other words, the uncertainties in the continuous dynamics are caused by an imperfect modeling of the drag.

The bounce happens when the state enters the guard set $\{y\leq 0,\dot{y}<0 \}$ in a deterministic fashion.
As the discrete jump is triggered spontaneously by a rate function in the GSHS formulation considered in this paper, we define the rate function to approximate the effects of the guard as follows \cite{hespanha2004stochastic}:
\begin{equation} \label{eqn:RateBall}
	\lambda(y,\dot{y}) = \begin{cases}
		100, &\text{if } y<0, \dot{y}<0 \\
		30, &\text{if } y=0, \dot{y}<0 \\
		0, &\text{otherwise}
	\end{cases}.
\end{equation}
In other words, (\ref{eqn:RateBall}) implies the jump never occurs when $y>0$ or $\dot{y}\geq0$, and it is very likely to occur when $y\leq 0$ and $\dot{y}<0$.
As the proposed computational method relies on a fixed time step, the height may have a negative value temporarily. 
In the above equation, the transition rate for $y<0$ is greater, as we wish to quickly remedy the occasion of negative height. 
This issue can be addressed by time-stepping in the numerical integration of the Fourier parameters in the future work.

During the discrete jump, the hybrid state is reset to $y^+=|y^-|$, $\dot{y}^+ = -c\dot{y}^- + w_c$, where $w_c\sim\mathcal{N}(0,\sigma_c^2)$ and $\mathcal{N}$ stands for the Gaussian distribution.
Thus, the reset kernel is given as
\begin{equation} \label{eqn:KernelBall}
	\kappa(y^-,\dot{y}^-,y^+,\dot{y}^+) = \frac{\delta(y^+-|y^-|)}{\sqrt{2\pi}\sigma_c}e^{-\frac{(\dot{y}^++c\dot{y}^-)^2}{2\sigma_c^2}},
\end{equation}
where $\delta$ is Dirac delta function, $c$ is the coefficient of restitution associated with the collision, and $s$, $s^{\prime}$ are omitted because there is only one discrete mode.
Therefore, the uncertainties in the discrete transitions correspond to the randomness in the velocity immediately after bouncing.
In (\ref{eqn:KernelBall}), $y^+$ is reset to $|y^-|$ because at the discretized time, if the ball goes below 0, the kernel will reset the position to its reflection $-y^-=|y^-|$ rather than 0, which is more realistic.

\begin{table}
	\centering
	\caption{Parameters of the bouncing ball model \label{table:ParaBall}}
	\begin{tabular}{|C{0.38in}|C{0.38in}|C{0.42in}|C{0.36in}|C{0.38in}|C{0.36in}|}
		\hline
		g & $\nu$ & $\sigma_v$ & $c$ & $\sigma_c$ & $\sigma_m$ \\
		\hline
		9.8$\,\frac{m}{s^2}$ & 0.05$\,\frac{1}{m}$ & 0.01$\,\frac{\sqrt{s}}{m}$ & 0.95 & 0.5$\,\frac{m}{s}$ & 0.3$\,m$ \\
		\hline
	\end{tabular}
\end{table}

The parameters of the bouncing ball model used in this simulation are given in Table \ref{table:ParaBall}.
The simulation lasts until $t=6$ seconds, and the time-step is $0.025$ seconds.
The PDF $p(t,y,\dot{y})$ is assumed to be supported on a rectangle $[-2.5,2.5]\times[-8,8]$ in the two space parameters, which is discretized as a 100$\times$100 grid.
We also compare the proposed method for propagating the PDF of a bouncing ball to a Monte Carlo simulation, where 1,000,000 samples are used to approximate the PDF by counting samples in the grid blocks.
The simulation is implemented in MATLAB R2019a, and the DFT (\ref{eqn:DFT}) and IDFT (\ref{eqn:IDFT}) are carried out using the built-in \texttt{fft}, \texttt{fft2} and \texttt{ifft2} functions respectively.

For numerical implementation, the PDF of bouncing ball has a sharp change at $\{y=0\}$.
Though the spectral method is good at dealing with large space variations \cite{tadmor2012review}, it still poses a challenge, since the vanishing rate of the Fourier coefficients is directly related to the smoothness of the PDF \cite{stein2011fourier}.
The lost information in the high order harmonics neglected in (\ref{eqn:FPContDFT}) could accumulate numerical errors with time.
Therefore, after each step, we set the PDF to exact 0 if $p(y,\dot{y})<3\times10^{-3}$ (including all negative values) and re-normalize it to eliminate the small ripples caused by the higher order terms omitted in (\ref{eqn:FPContDFT}).

\subsection{Uncertainty Propagation Results}

We first propagate the PDF assuming that there is no measurement available. 
The initial distributions of the position and velocity for uncertainty propagation are set to be independent, as $y\sim\mathcal{N}(1.5,0.2^2)$, and $\dot{y}\sim\mathcal{N}(0,0.5^2)$.
In other words, the ball is dropped from $1.5\,\mathrm{m}$ with zero initial speed with some uncertainties.

\begin{table}
	\centering
	\caption{Simulation time of the bouncing ball using two methods. \label{table:timeBallProp}}
	\begin{tabular}{|L{1.1in}|C{0.85in}|C{0.85in}|}
		\hline
		& total time & step time \\
		\hline
		Spectral-splitting & 2$\,$min 38$\,$s & 0.076$\,$s \\
		\hline
		Monte Carlo & 10$\,$min 44$\,$s & 2.67$\,$s \\
		\hline
	\end{tabular}
	\raggedright
	\textit{Note:} step time indicates the CPU time spent to propagate the PDF in one time-step.
\end{table}

\begin{figure}
	\centering
	\includegraphics{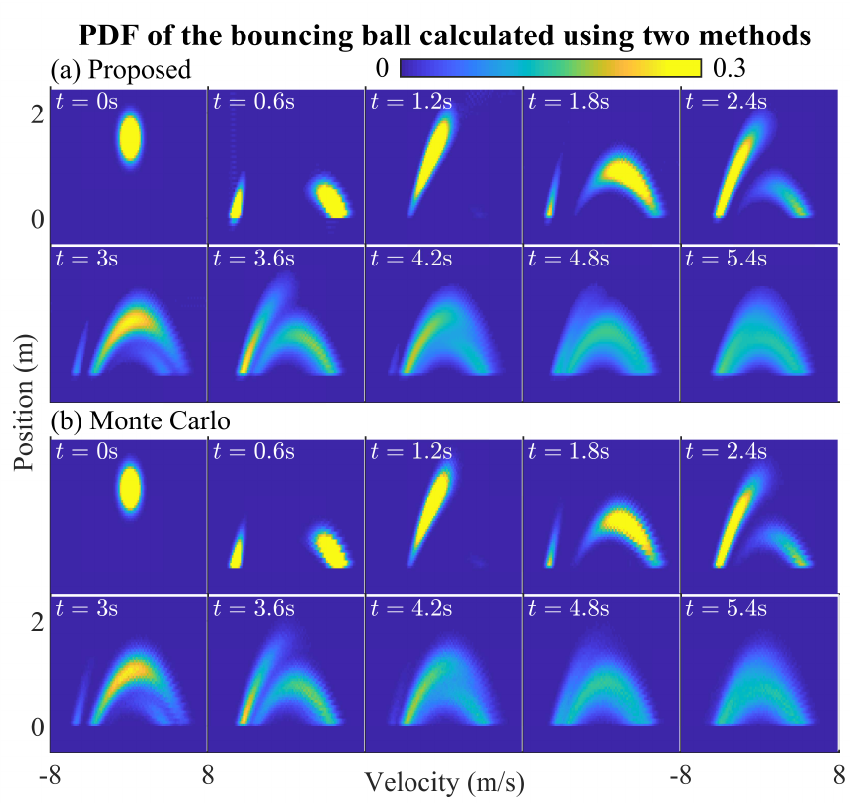}
	\caption{The evolution of the PDF of the bouncing ball model calculated by (a) the proposed spectral-splitting method, and (b) Monte Carlo method from 0s to 5.4s. \label{fig:ballPDF}}
\end{figure}

The PDF of the bouncing ball model calculated from the proposed method is compared with the Monte Carlo method in Fig. \ref{fig:ballPDF}.
At $t=0\,$s, the ball is at about $1.5\,\mathrm{m}$ with the velocity of approximately $0\,\mathrm{m/s}$.
Then, the height reduces and the velocity decreases to be negative due to the gravity, and at about $t=0.6\,$s, the first bounce occurs and its velocity jumps to be positive (from left to right in the figure).
Interestingly, the non-zero values of the probability density are split into two parts: the left group represents the cases before any bouncing, and the right group after the first bouncing. 
Due to the additional uncertainties in the discrete transition, the latter is more dispersed. 
Then the position begins to rise and the velocity is decelerated by gravity.
These are repeated after the velocity becomes reduced to zero again.
Due to the accumulation of uncertainties, the PDF gradually becomes more and more dispersed as time increases.
At $t=6\,$s, the ball has roughly bounced 6-7 times, and the distribution becomes a shape of the letter ``U" inverted.
Such a distribution cannot be approximated by the common Gaussian distribution.

Furthermore, the PDF calculated from the spectral-splitting method is almost identical to the Monte Carlo simulation, thereby verifies its accuracy.
The computation time of the two methods is compared in Table \ref{table:timeBallProp}.
Because the burdensome calculation of the matrix exponential in \eqref{eqn:ContSol} and \eqref{eqn:DiscSol} can be reused in each time step, the splitting-spectral method is much more efficient in propagating the density over one time-step than the Monte Carlo method, which needs to propagate a lot of samples and count them in each grid block.

\subsection{Estimation Results}

For the estimation simulation, the initial distribution is set to be uniform on $[0,2.5]\times[-8,8]$, implying there is completely no prior knowledge about the state.
The height of the ball $y$ is assumed to be measured by a sensor, and the measurement function is
\begin{equation}
	z = y+v,
\end{equation}
where $v\sim\mathcal{N}(0,\sigma_m^2)$ is the measurement noise.
Therefore, the likelihood function is given as
\begin{equation}
	p\left(z \lvert y,\dot{y}\right) = \frac{1}{\sqrt{2\pi}\sigma_m}e^{-\frac{(z-y)^2}{2\sigma_m^2}}.
\end{equation}
Measurements are assumed to be available in every time step.
Although the velocity $\dot{y}$ is not directly measured, it builds up correlations with the height $y$ through the dynamics, so the measurement can still make proper corrections to $\dot{y}$.

Because multiplying the likelihood function could potentially magnify the numerical errors caused by neglecting the higher order harmonics, we adopt a stronger criterion than in the uncertainty propagation simulation, i.e. all densities smaller than 1/40 of the peak density are set to zero, to eliminate the ripples in the prior density.
In addition, the posterior distribution becomes multi-modal around collisions, therefore we choose the state maximizing the PDF as the estimated state.

The proposed method is compared with an ISR particle filter with 1,000,000 particles.
The resampling step is implemented in every time step by the systematic resampling scheme \cite{arulampalam2002tutorial}.
In each step, the PDF is recovered by counting particles over the grid.
Sixty simulations with randomly generated true sample paths and noises are used to test the two estimation methods.

\begin{figure}
	\centering
	\includegraphics{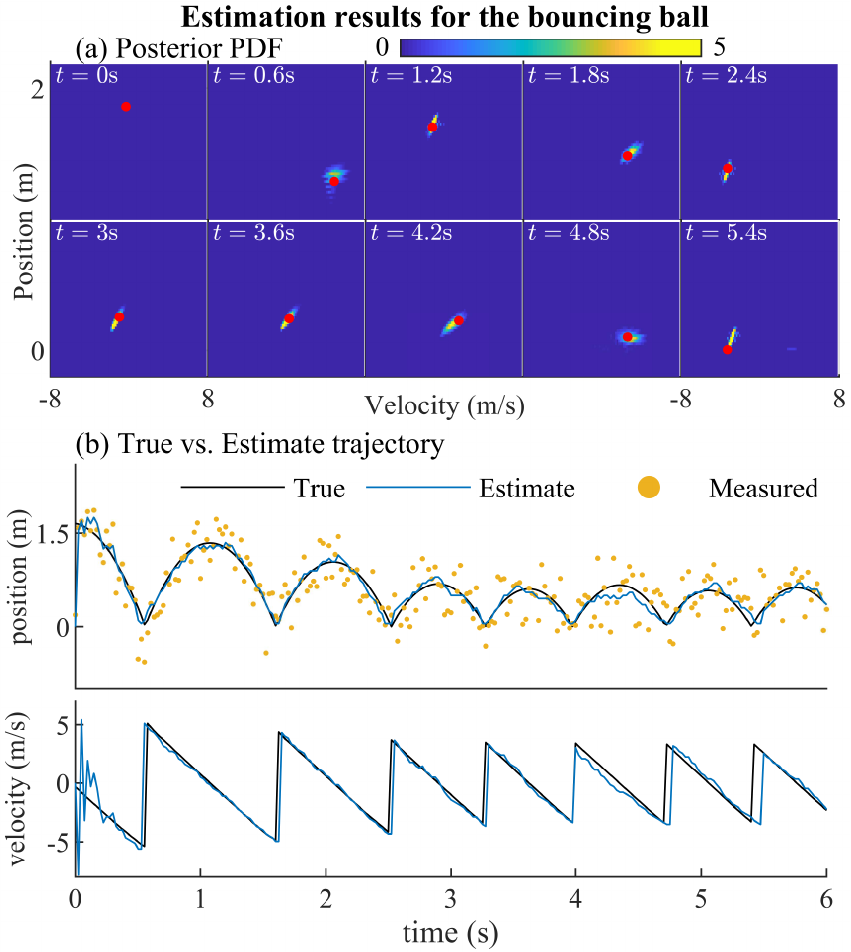}
	\caption{Estimation results of the bouncing ball. (a) The posterior PDF from 0s to 5.4s, the red dots represent the true state. (b) The true, estimated, measured position and velocity of the ball. \label{fig:ballEst}}
\end{figure}

The results of an example simulation using the proposed splitting-spectral method are shown in Fig. \ref{fig:ballEst}.
In subfigure (a), the posterior PDF is presented, which converges quickly from uniform distribution at $t=0\,$s to a concentrated distribution peaked around the true state.
In subfigure (b), the true, estimated and measured hybrid states are compared.
Again, the estimated position and velocity quickly converge around the true state, though the velocity is not directly measured.
Furthermore, even though the initial state is completely unknown, the estimated position is much more accurate than the measurement, which demonstrates the contributions from the model dynamics.

\begin{table}
	\centering
	\caption{Comparison between the proposed and particle filters \label{tab:SSvsPF}}
	\begin{tabular}{|c|c|c|c|}
		\hline
		 & \thead{position \\ error (m)} & \thead{velocity \\ error (m/s)} & \thead{step \\ time (s)} \\ \hline
		Proposed & 0.091$\pm$0.017 & 0.68$\pm$0.14 & 0.098$\pm$0.029 \\ \hline
		Particle & \thead{0.091$\pm$0.015 \\ $\scriptstyle (p=0.41)$} & \thead{0.73$\pm$0.14 \\ $\scriptstyle (p<0.001)$} & \thead{2.99$\pm$0.05 \\ $\scriptstyle (p<0.001)$} \\ \hline 
	\end{tabular}
\end{table}

The comparisons of estimation accuracy and step time with the particle filter are presented in Table \ref{tab:SSvsPF}.
The errors and computation time are averaged across time steps in each simulation, and further averaged across sixty simulations.
The $p$-values are obtained from paired $t$-tests ($N=60,\,\alpha=0.05$) to indicate statistical differences.
It is shown that the spectral-splitting method has more accurate velocity estimates and is more computationally efficient than this particular particle filter.
The major reason for the particle filter to be less accurate in velocity is during the converging process, the uncertainty is large, therefore some particles are outside of the grid but counted as in the boundary blocks to construct the PDF, making the PDF peak at the boundary.


\section{Application to Dubins Vehicle} \label{section:car}

In the bouncing ball example, there is only one discrete mode.
To demonstrate the efficacy of the proposed method, here we apply it to a Dubins vehicle model with multiple modes \cite{dubins1957curves}.
Dubins vehicle is a simple kinematic model that moves along its heading direction at  a constant speed with a turning control, and 
it is commonly used for motion planning of wheeled robots in robotics or airplanes in air traffic control. 
We consider three discrete modes: moving forward, turning left and turning right.
However, the continuous state is not reset during the discrete jump, and consequently, the Dubins vehicle model is a switching diffusion process.

\subsection{Dubins Vehicle Model}

The Dubins vehicle is assumed to be moving on a horizontal plane with fixed speed and a controllable heading angle.
The continuous state is $r=(y,\theta)$, where $y=(y_1,y_2)\in\mathbb{R}^2$ is the horizontal position of the vehicle, and $\theta\in\mathbb{R}/2\pi$ is the heading angle.
The continuous dynamics is given as the following SDE:
\begin{equation}
	\text{d}\left[\begin{matrix}y_1\\y_2\\\theta\end{matrix}\right] = \left[\begin{matrix}v\cos\theta\\v\sin\theta\\u(s)\end{matrix}\right]\text{d}t + \left[\begin{matrix}0\\0\\\sigma_u\end{matrix}\right]\text{d}W_t,
\end{equation}
where $v\in\mathbb{R}^+$ is the constant speed, $u(s)$ is a set of turning rates disturbed by a Gaussian white noise with strength $\sigma_u$.
The turning rate $u(s)$ is associated with three discrete modes $S=\{1,2,3\}$, corresponding to moving forward, turning left, and turning right, as follows:
\begin{align}
    u(s) = \begin{cases}
        0 & s=1 \\
        a & s=2, \\
        -a & s=3
    \end{cases}
\end{align}
where $a\in\mathbb{R}^+$ is a constant value.
The vehicle cannot make a sharp turn, but follows an arc with radius $v/a$ during turning.
It should be noted the heading angle $\theta$ belongs to a compact circular space homomorphic to $\mathbb{S}^1$, where the Fourier series is naturally defined.

\begin{figure}
	\centering
	\includegraphics{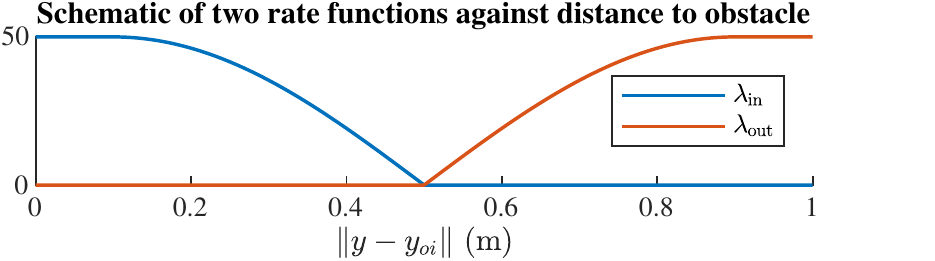}
	\caption{The ``in" and ``out" rate functions plotted against the distance to the obstacle. \label{fig:carRate}}
\end{figure}

\begin{table}
	\centering
	\caption{Parameters of the Dubins vehicle model \label{table:ParaCar}}
	\begin{tabular}{|C{0.7in}|C{0.6in}|C{0.65in}|C{0.65in}|}
		\hline
		$v$ & $a$ & $\sigma_u$ & $y_{oi,1}\,$ \\
		\hline
		1 $\,$ m/s & 2$\,$ rad/s & 0.2$\,$ rad/$\sqrt{s}$ & 0,$\,$1,$\,$1$\,$m \\
		\hhline{|=|=|=|=|}
		$y_{oi,2}\,$ & $d$ & $\sigma_l$ & $\kappa_l$ \\
		\hline
		0,$\,$-1.5,$\,$1.5$\,$m & 0.5$\,$m & 0.5$\,$m & 30 \\
		\hline
	\end{tabular}
\end{table}

We place three obstacles on the plane at $\{y_{oi}\}_{i=1}^3$, and
the discrete jump is designed to avoid collisions with these obstacles as follows.
When the vehicle comes sufficiently close to any obstacle, and enters the set $\{\ensuremath{\left\|y-y_{oi}\right\|}<d\}$ for $d>0$, it turns either left or right depending on whether the obstacle is on its right side of the path or on the left side. 
After it moves sufficiently far away from the obstacle and leaves the set, it stops turning and continues to move forward.
More specifically, we use the following two rate functions with non-intersecting supports to model the two discrete transitions.
\begin{equation*}
	\lambda_{\text{in}}(y) = \begin{cases}
		50, \hspace{0.4in} \text{if} \; \exists i, \; \ensuremath{\left\|y-y_{oi}\right\|}<0.1 \\
		50\sin\left(\frac{0.5-\ensuremath{\left\|y-y_{oi}\right\|}}{0.4} \frac{\pi}{2}\right), \\
		\hspace{0.6in} \text{if} \; \exists i, \; 0.1\leqslant\ensuremath{\left\|y-y_{oi}\right\|}<0.5 \\
		0, \hspace{0.47in} \text{if} \; \forall i, \; \ensuremath{\left\|y-y_{oi}\right\|}\geqslant 0.5
	\end{cases},
\end{equation*}
\begin{equation*}
	\lambda_{\text{out}}(y) = \begin{cases}
		50, \hspace{0.4in} \text{if} \; \forall i, \; \ensuremath{\left\|y-y_{oi}\right\|}>0.9 \\
		50\sin\left(\frac{\ensuremath{\left\|y-y_{oi}\right\|}-0.5}{0.4} \frac{\pi}{2}\right), \\
		\hspace{0.6in} \text{if} \; \exists i, \; 0.5<\ensuremath{\left\|y-y_{oi}\right\|}\leqslant0.9 \\
		0, \hspace{0.47in} \text{if} \; \exists i, \; \ensuremath{\left\|y-y_{oi}\right\|}\leqslant 0.5
	\end{cases},
\end{equation*}
and
\begin{equation}
	\lambda(y,s) = \begin{cases}
		\lambda_{\text{in}}(y), &\text{if} \; s=1 \\
		\lambda_{\text{out}}(y), &\text{if} \; s=2 \; \text{or} \; 3
	\end{cases}.
\end{equation}
The non-intersecting supports guarantee the discrete mode does not switch back and forth between going straight and turning at the edge of the obstacle.
The \texttt{sin} function is used to connect the two rates (0 and 50) as depicted in Fig. \ref{fig:carRate}.

Let $Y_i=\left\{\ensuremath{\left\|y-y_{oi}\right\|}<d\right\}$ denote the set where the vehicle should make a turn, and $\overline{Y}=\mathbb{R}^2\setminus\cup_{i=1}^3Y_i$ denote the set where the vehicle should move forward.
Because the discrete jump does not reset the continuous state, the kernel function is expressed as the product of Dirac delta functions:
\begin{equation*}
	\tilde{\kappa}(r,s^-,s^+) = \begin{cases}
		\delta(s^+-1), &\text{if} \; s^-=2 \; \text{or} \; 3 \; \text{and} \; y\in\overline{Y}\\
		\delta(s^+-2), &\parbox[t]{1.4in}{if$ \, s^-=1, \, \exists i, \, 
			y\in Y_i, \; $ and $ \; -\pi\leqslant\theta_{oi}-\theta<0 $} \\
		\delta(s^+-3), &\parbox[t]{1.4in}{if$\; s^-=1, \; \exists i, \; 
			y\in Y_i, \; $ and $ \; 0\leqslant\theta_{oi}-\theta<\pi $}
	\end{cases},
\end{equation*}
and
\begin{equation} \label{eqn:KernelCar}
	\kappa(r^-,s^-,r^+,s^+) = \delta(r^+-r^-)\tilde{\kappa}(r^-,s^-,s^+),
\end{equation}
where $\theta_{oi}=\texttt{atan2}(y_{oi,2}-y_2,y_{oi,1}-y_1)$ corresponds to the direction toward the $i$-th obstacle from the vehicle.
As noted in Section \ref{section:quadrature}, the discrete dynamics can be propagated using the simplified ODE system (\ref{eqn:FPDiscNumSDP}).

\begin{figure}
	\centering
	\includegraphics{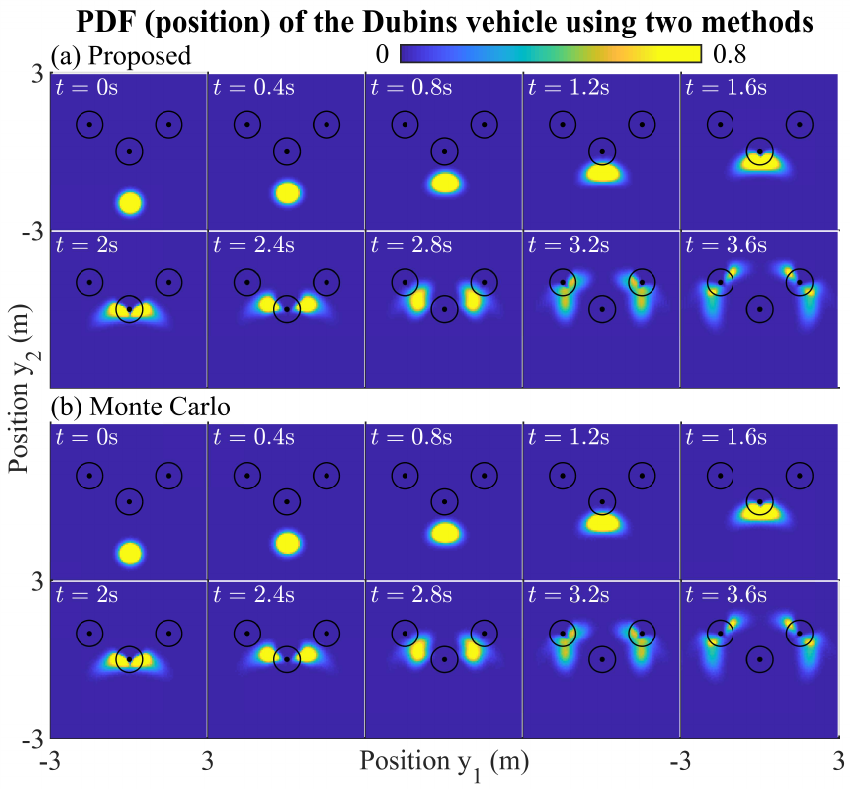}
	\caption{The evolution of the position marginal PDF of the Dubins vehicle model calculated by (a) the proposed spectral-splitting method, and (b) Monte Carlo method from 0s to 3.6s. The black dots and circles are the obstacles and their affecting ranges. \label{fig:carPDFPosition}}
\end{figure}

The parameters of the Dubins vehicle model used in this simulation are given in Table \ref{table:ParaCar}.
The simulation lasts until $t=4$ seconds, and the time-step is $0.025$ seconds. 
The PDF $p(t,y_1,y_2,\theta,s)$ is supposed to be supported on $[-3,3]\times[-3,3]\times[0,2\pi)$ in the continuous state space, which is discretized as a $100\times100\times50$ grid.
We also compared the proposed uncertainty propagation algorithm to a Monte Carlo simulation with 1,000,000 samples.

\subsection{Uncertainty Propagation Results}

\begin{figure}
	\centering
	\includegraphics{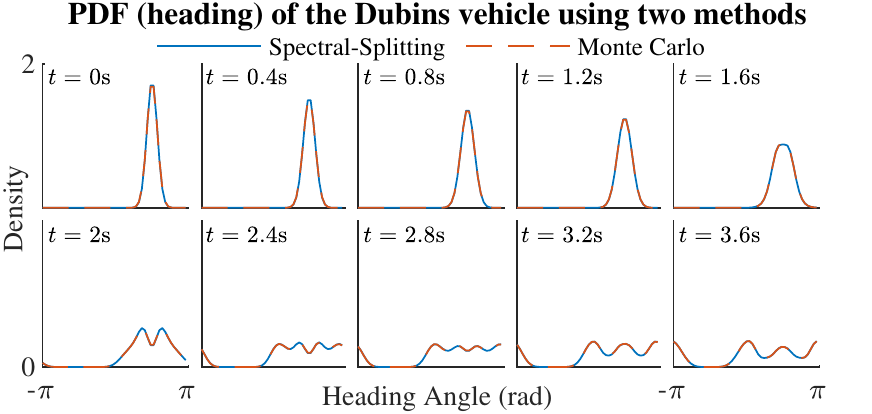}
	\caption{Heading marginal PDF of the Dubins vehicle model. \label{fig:carPDFHeading}}
\end{figure}

\begin{figure}
	\centering
	\includegraphics{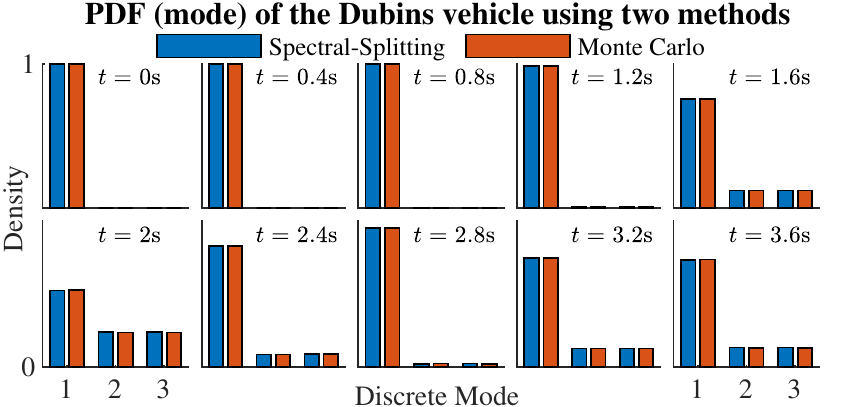}
	\caption{Discrete mode marginal PDF of the Dubins vehicle model. \label{fig:carPDFMode}}
\end{figure}

We first propagate the PDF without any measurements.
The initial distribution is set as $y_1\sim\mathcal{N}(0,0.2^2)$, $y_2\sim\mathcal{N}(-2,0.2^2)$, $\theta\sim \mathrm{VM}(\pi/2,20)$, $\mathbb{P}(s=1)=1$, and they are independent, where VM$(\mu,\kappa)$ stands for the von Mises distribution.
The von Mises distribution is the analog of Gaussian distribution defined in the compact circular space $\mathbb{R}/2\pi$ \cite{mardia2009directional}, and has density function:
\begin{equation}
	p(\theta;\mu,\kappa) = \frac{1}{2\pi I_0(\kappa)}e^{\kappa\cos(\theta-\mu)},
\end{equation}
where $I_0$ is the modified Bessel function of the first kind
In other words, the vehicle starts at around $y=(0,-2)$, with initial heading angle at around $\pi/2$, and in discrete mode $s=1$ (moving forward).

The propagated marginal PDF of the Dubins vehicle model on position, heading angle and discrete mode are shown in Fig. \ref{fig:carPDFPosition}-\ref{fig:carPDFMode} respectively.
During $t=0\,\mathrm{s}$ to $t=1.2\,\mathrm{s}$, the vehicle moves forward from its initial position. The heading angle becomes more dispersed as the uncertainty in turning rate accumulates, and the position density becomes arc-shaped.
At about $t=1.6\,\mathrm{s}$, the vehicle encounters the first obstacle, and the density bifurcates into two directions: one half part turns left ($s=2$) and the other half turns right ($s=3$).
Then after the vehicle leaves the first obstacle at around $t=2.8\,\mathrm{s}$, it continues to move forward ($s=1$).
Finally at $t=3.6\,\mathrm{s}$, the vehicle meets another obstacle, and the two branches of density bifurcate again into four branches.

Furthermore, the PDF calculated using the spectral-splitting approach matches the Monte Carlo simulation almost exactly.
The computation time of the two methods are compared in Table \ref{table:timeCarProp}.
It is shown that even in this more complex model, the proposed method is still more efficient than the Monte Carlo method with the same number of samples as in the bouncing ball example.

\begin{table}
	\centering
	\caption{Simulation time of the Dubins vehicle using two methods. \label{table:timeCarProp}}
	\begin{tabular}{|L{1in}|C{0.9in}|C{0.9in}|}
		\hline
		& total time & step time \\
		\hline
		Spectral-splitting & 10$\,$min 30$\,$s & 3.07$\,$s \\
		\hline
		Monte Carlo & 20$\,$min 39$\,$s & 7.70$\,$s \\
		\hline
	\end{tabular}
\end{table}

\subsection{Estimation Results}

For Bayesian estimation, it is assumed that the position of the vehicle is measured by a Lidar located at $y_l=(0,-3)$.
The Lidar measures the distance $d_l\in\mathbb{R}^+$ from itself to the vehicle, and the angle $\theta_l\in\mathbb{R}/2\pi$ formed by the vehicle, itself and positive $y_1$ axis.
Therefore, the measurement function is given as
\begin{equation}
	z = \left[\begin{matrix} d_l \\ \theta_l \end{matrix}\right] = 
	\left[\begin{matrix}
		\ensuremath{\left\|y-y_l\right\|} + v_d \\
		\mathrm{atan2}(y_2-y_{l2},y_1-y_{l1}) + v_\theta
	\end{matrix}\right],
\end{equation}
where $v_d\sim\mathcal{N}(0,\sigma_l^2)$, and $v_\theta\sim \mathrm{VM}(0,\kappa_l)$ are measurement noises.
Therefore, the likelihood function is
\begin{equation}
	p(d_l,\theta_l\lvert y,\theta,s) = \frac{1}{(2\pi)^{3/2}\sigma_lI_0(\kappa_l)} e^{-\frac{(d_l-\tilde{d}_l)^2}{2\sigma_l^2}+\kappa_l\cos(\theta_l-\tilde{\theta}_l)},
\end{equation}
where $\tilde{d}_l=\ensuremath{\left\|y-y_l\right\|}$, $\tilde{\theta}_l=\mathrm{atan2}(y_2-y_{l2},y_1-y_{l1})$.
The measurement is assumed to come at every time step.
It should be noted the heading angle and discrete mode are unobserved states, but they build correlations with the position through the dynamics equations, thus they can be corrected by the measurement as well.

\begin{figure}
	\centering
	\includegraphics{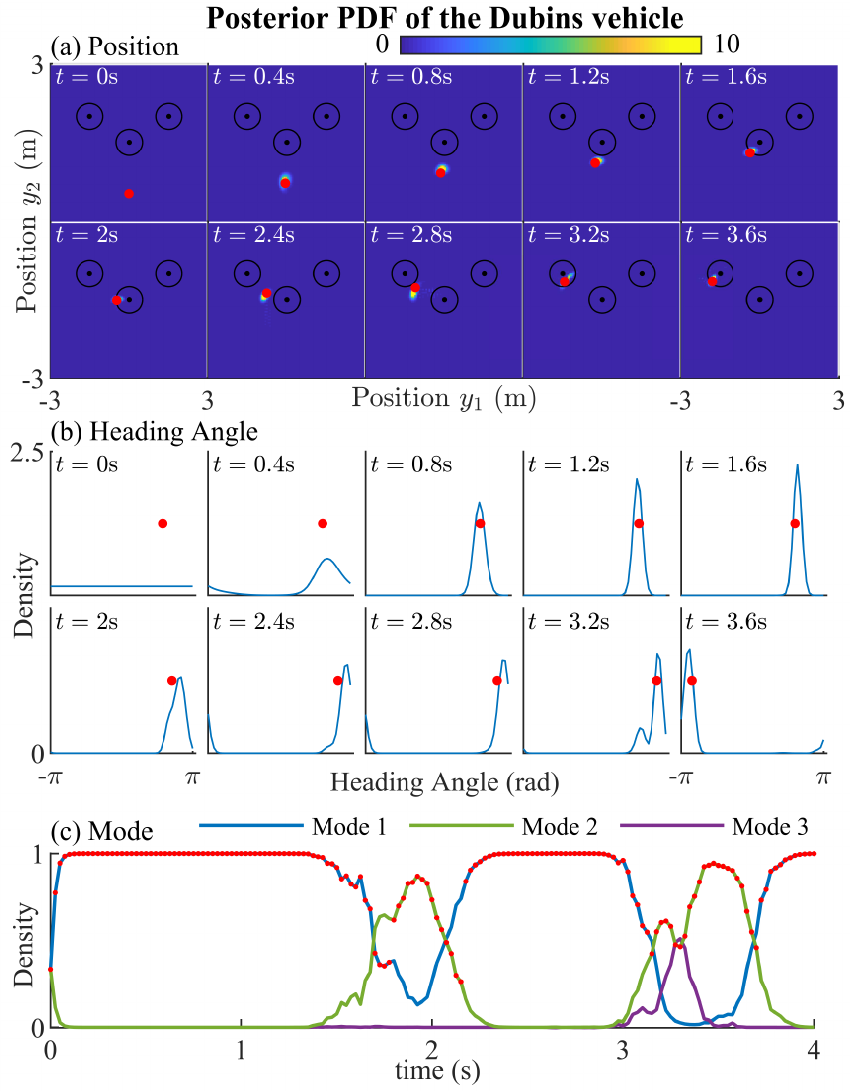}
	\caption{Marginal posterior PDF of the Dubins vehicle with Lidar measurements on: (a) position, (b) heading angle, and (c) discrete mode. Red dots represent the true state. \label{fig:carEstPDF}}
\end{figure}

\begin{figure}
	\centering
	\includegraphics{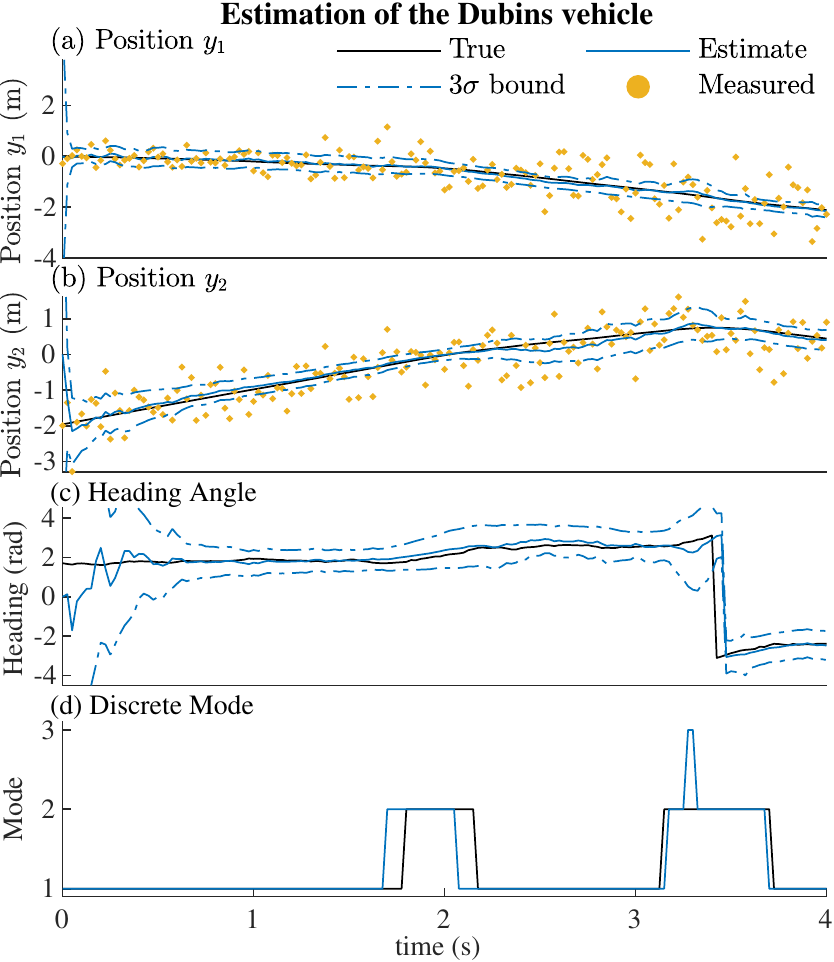}
	\caption{Estimation of the hybrid state of the Dubins vehicle: (a) position $y_1$, (b) position $y_2$, (c) heading angle $\theta$, and (d) discrete mode $s$. The dash-dot lines represent the triple standard deviation ($\sigma$) bound calculated from the marginal PDF, the $\sigma$ in (c) is computed as the circular standard deviation, see \cite{mardia2009directional} p.30. \label{fig:carEstTraj}}
\end{figure}

The initial distribution for the estimation simulation is assumed to be uniform on the continuous grid and in discrete modes, i.e., there is no prior information about the vehicle's initial state.
We use the expected value as the position estimate, the mean direction (see \cite{mardia2009directional} p.29) as the heading estimate, and the most probable mode as the discrete state estimate.

The proposed spectral-splitting method is compared with an ISR particle filter with 1,000,000 particles and an IMM approach modified to handle the continuous state dependent transition.
The resampling for the particle filter is implemented in every time step by the systematic method.
For the IMM approach, the discrete state transition matrix is obtained by numerically integration given by the following equation
\begin{align} \label{eqn:IMMTran}
p(s_k&=j|s_{k-1}=i) = \nonumber \\
&\int_{r\in\mathbb{R}^3}(1-e^{-\lambda(y,i)\Delta t})\tilde{\kappa}(r,i,j)p_{k-1}(r,i)\mathrm{d}r,
\end{align}
where $p_{k-1}(r,i)$ is the Gaussian density in the Kalman filter of discrete mode $i$.
To match the proposed and particle filters, the order of continuous and discrete propagation in IMM is reversed compared to the original algorithm \cite{blom1988interacting}, i.e. the continuous propagation is applied first, which we assume is a minor modification to facilitate comparisons.
Sixty randomly generated true sample paths and noises are used to test these algorithms.

An example simulation of the Dubins vehicle estimated by the spectral-splitting method is shown in Fig. \ref{fig:carEstPDF} and \ref{fig:carEstTraj}.
Fig. \ref{fig:carEstPDF} presents the posterior marginal density, it can be seen
the position PDF quickly converges from uniform at $t=0\,\mathrm{s}$ to a concentrated density peaked around the true position.
Because the heading angle is not directly observed, it converges slower.
The discrete mode also converges from equally distributed point masses to exclusively at $s=1$ (moving forward) in a short time.

The true, estimated and measured states are compared in Fig. \ref{fig:carEstTraj}, with 3$\sigma$ bounds shown for continuous states.
After convergence, the position and heading estimates follow closely to the true trajectories, even when the vehicle moves far away from the Lidar and the noise in angle measurement gradually degrades the position accuracy (Fig. \ref{fig:carEstTraj}(a)).
Upon transitions between going forward and turning, due to the noisy measurement, density in the false branch during bifurcation sometimes becomes larger, causing false estimates of the discrete mode, and increased uncertainties in continuous states (around 3.3$s$).

The comparisons in estimation accuracy and computation time of the three filters are shown in Table \ref{tab:carCom}.
The error terms and computation time are similarly defined as in the bouncing ball example, except for the mode error which represents the percentage of false estimates in a simulation averaged across sixty simulations.
The spectral-splitting method is significantly better than the IMM in both estimation accuracy and computation efficiency.
The particle filter achieves similar estimation accuracy except for the mode estimate, which is slightly but significantly worse than the proposed method.
Moreover, the computation time of the proposed method is less than this particular particle filter.

\begin{table}
	\centering
    \caption{Comparison with~\cite{liu2014hybrid}: accuracy and computation time  \label{table:chasing}}

	\begin{tabular}{|c|c|c|c|}
		\hline
		\thead{grid \\ size (m)} & \thead{distance \\ error (m)} & \thead{mode \\ accuracy (\%)} & \thead{step \\ time (ms)} \\ \hline
		5 & 18.3$\pm$3.8 & 85$\pm$4 & 110 \\ \hline
		10 & 18.3$\pm$3.8 & 85$\pm$4 & 28 \\ \hline
		20 & 18.6$\pm$3.8 & 85$\pm$4 & 6.9 \\ \hline
		25 & 19.4$\pm$3.8 & 85$\pm$4 & 4.6 \\ \hline
		30 & 20.9$\pm$3.8 & 85$\pm$4 & 3.4 \\ \hline
		40 & 26.3$\pm$4.1 & 85$\pm$4 & 2.0 \\ \hline
	\end{tabular}
\end{table}


\begin{table*}
	\centering
	\caption{Accuracy and computation time comparisons of spectral-splitting, particle and IMM filters. \label{tab:carCom}}
	\begin{tabular}{|c|c|c|c|c|c|}
		\hline
		& $y_1$ error (m) & $y_2$ error (m) & $\theta$ error (rad) & mode error (\%) & computation time (s) \\ \hline
		Proposed & 0.083$\pm$0.025 & 0.092$\pm$0.024 & 0.35$\pm$0.10 & 8.3$\pm$6.1 & 3.8$\pm$0.1 \\ \hline
		Particle & 0.083$\pm$0.024 $\scriptstyle (p=0.67)$ & 0.093$\pm$0.024 $\scriptstyle (p=0.06)$ & 0.35$\pm$0.10 $\scriptstyle (p=0.08)$ & 8.6$\pm$6.2 $\scriptstyle (p=0.02)$ & 8.9$\pm$0.1 \\ \hline
		IMM & 0.156$\pm$0.086 & 0.194$\pm$0.15 & 0.70$\pm$0.40 & 13.5$\pm$8.6 & 10.5$\pm$0.2 \\ \hline
	\end{tabular}
	\raggedright
	The $p$-values between the proposed method and the other two methods are below 0.001 except those indicated.
\end{table*}

\section{Comparison with \cite{liu2014hybrid}}
\begin{figure}
	\includegraphics{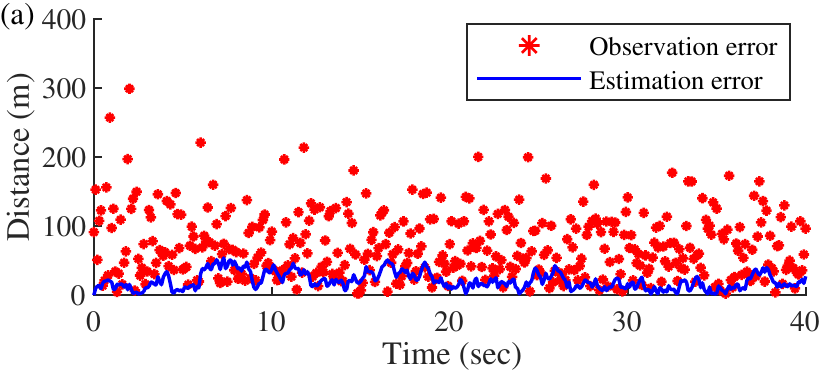}
	\includegraphics{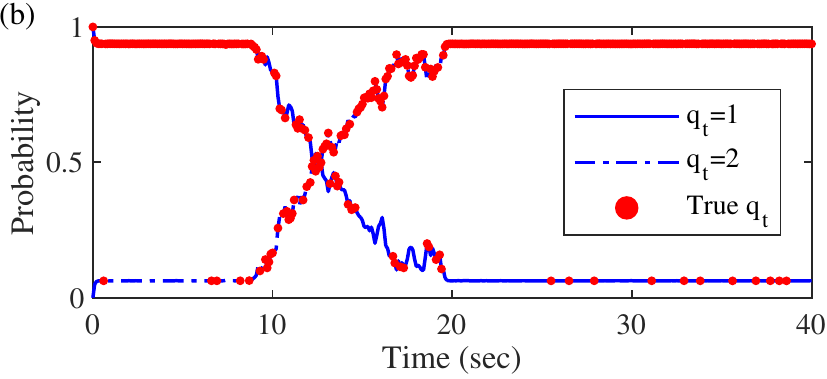}
	\caption{Comparison with \cite{liu2014hybrid}: (a) Estimation vs. observation error of the continuous state. (b) Estimated discrete state probability vs. true discrete state. \label{fig:TwoVehicle}}
\end{figure}

Since this paper is similar to \cite{liu2014hybrid} in the idea of solving the FP equation \eqref{eqn:FP}, but differs in numerical method (spectral vs. finite difference), in this section we apply the proposed spectral-splitting method to the aircraft chasing example in \cite{liu2014hybrid} (Section V) to conduct a comparison study.
The model and parameters are chosen exactly the same as \cite{liu2014hybrid}, with the grid size ranging from 5 meters to 40 meters.
Sixty Monte Carlo simulations were conducted, the accuracy and computation step time are shown in Table \ref{table:chasing}.
Readers are invited to compare these results with Fig. 9 and Table II in \cite{liu2014hybrid}.
The continuous state error and estimated discrete state probability from one Monte Carlo simulation with 5m grid size are shown in Fig. \ref{fig:TwoVehicle}.

In short, with the same grid size, the estimation of continuous state using the proposed method is more accurate than the method in \cite{liu2014hybrid}, at the cost of longer computation time.
This is anticipated because as suggested in \cite{tadmor2012review}, the spectral method is typically more accurate but computational intensive than the finite difference method.
The reason is that the spectral method is ``global'', i.e. it uses all the information from grid to approximate the differentiation; whereas the finite difference method only uses the nearby grid points.
The accuracy of discrete state estimation seems to be worse than the method in \cite{liu2014hybrid}, however, this is more possibly due to different implementations of generating the true trajectory, as there are a lot of oscillations between the two discrete states in the true trajectory we generated (see Fig. \ref{fig:TwoVehicle}(b)) compared to that in \cite{liu2014hybrid}.

\section{Discussion} \label{sec:discussion}
The proposed spectral-splitting method and particle filter are the so called ``exact'' Bayesian filtering approaches, which approximate the PDF directly instead of using parametric assumptions, such as mean and variance for Gaussian distributions.
Even though the noises are assumed to be Gaussian in both examples, the solutions to the FP equations are no longer Gaussian, and sometimes become even multi-modal, because of discrete transitions and non-linearities involved in the continuous SDEs.
Therefore, it is expected that the IMM approach which relies on the Gaussian assumption and linearization is much less accurate than the other two filters in the Dubins vehicle example.
Furthermore, due to the numerical integration in \eqref{eqn:IMMTran} for IMM, it does not show any computational advantage over the proposed method.
Approaches have been proposed to avoid this numerical integration \cite{seah2009state}, but at the cost of even compromised estimation accuracy.

The sequential importance resampling (SIR) particle filter is capable of estimating all kinds of GSHS, and does not assume any specific distribution or suffer linearization errors.
The performance of a particle filter is highly dependent on implementation, especially the number of particles and the frequency of resampling.
The particle filters tested in this paper have a million particles and perform resampling every step, which seems a little bit redundant, especially for the simple bouncing ball model.
Nevertheless, one of the motivations of the proposed method is to calculate the PDF of hybrid states, which encodes much more stochastic information than mere the estimates.
Even though the PDF can also be recovered from particles, this operation is computational intensive, which makes the particle filter less efficient.
More importantly, when the distribution is widely dispersed, the particles are scattered into a lot of, instead of a few, grid blocks.
This makes the number of particles needed to give a reliable density function increase drastically compared to a concentrated distribution.
In the simulation, it is observed that even with a million particles, the heading density of the Dubins vehicle model is very rough in the first a few steps, because the distribution is converging from a uniform density and is very widely dispersed.

The proposed spectral-splitting method originates from numerically solving the Fokker-Planck equation \eqref{eqn:FP}.
Similar to other numerical methods for PDEs, including the finite difference method adopted in \cite{liu2014hybrid}, the guard set presents significant challenges, as it makes the PDF in continuous state space discontinuous, where the approximation of derivatives breaks down.
Compared to the finite difference method, the spectral method is regarded as ``global'', since it uses all the information from the grid to approximate the differentiation, whereas the finite difference method only uses the grid points nearby.
Therefore, the spectral method adopted in this paper is typically more accurate with the same grid size but at the cost of higher computational demand.
Also, the proposed method suffers the so called curse of dimensionality, as it is clearly seen that the grid size increases exponentially with the dimension of the continuous space.
On the other hand, the computation complexity only grows linearly with the number of discrete states.
Therefore, the proposed method is more suitable for a GSHS with a large number of discrete modes but a few continuous dimensions.

One potential way to mitigate the two aforementioned drawbacks is to further split the continuous state space into partially overlapped tiles.
Since in the estimation problem, the PDF is usually supported in a small region, it is highly possible to be supported in one of the tiles, and the FP equation can be only solved in this tile where a finer and smaller grid can be used.
This could avoid the ``waste of grids" when the PDF on most grid points are zero, and therefore carry little information of the distribution as seen in Fig. \ref{fig:ballEst} and \ref{fig:carEstPDF}.
Besides, by utilizing non-commutative harmonic analysis, the proposed method can be naturally extended to a compact Lie group as the continuous state space, such as the circular space and especially SO(3), SE(3) which are ubiquitous in robotics and aerospace applications.

\section{Conclusions} \label{section:conclusion}
In this paper, we propose a computational approach to conduct Bayesian estimation of a GSHS, consisting of an uncertainty propagation step and a correction step.
For the uncertainty propagation, we use spectral method to solve the partial-differential part of the FP equation, numerical integration to solve the integral part, and splitting method to combine the two solutions into the solution to the entire FP equation.
The propagated uncertainty, i.e. the PDF of GSHS, is used with a likelihood function in the Bayes' formula to estimate the hybrid state in an exact manner.
We test the proposed spectral-splitting method in two examples: the bouncing ball model and the Dubins vehicle model to demonstrate its effectiveness in different types of GSHS.
The PDF propagated by the proposed method is consistent with that propagated by a Monte Carlo simulation.
When measurements are available, the proposed method achieves similar estimation accuracy and higher efficiency compared with the particular particle filter tested, and is much better in both aspects than the Gaussian based IMM approach in the Dubins vehicle model.

\begin{ack}
This work has been supported in part by AFOSR under the grant FA9550-18-1-0288.
\end{ack}

\bibliographystyle{plain}
\bibliography{AUTHybrid}

\appendix
\section{Multi-dimensional continuous space} \label{appendix}
We give the multi-dimensional version of (\ref{eqn:FPContDFT}) using the notations from the last paragraph in Section \ref{section:Spectral}.
The proof is similar to (\ref{eqn:FPContDFT}) using multi-dimensional versions of differentiation and product formulae.
As it is straightforward and tedious, we skip the detailed proof.
But, the resulting ODE is presented at \eqref{eqn:FPContDFTMulti}.

\begin{figure*}
    \begin{align} \label{eqn:FPContDFTMulti}
        \frac{\mathrm{d}}{\mathrm{d}t}\hat{f}_{n_1,\ldots,n_{N_r}}[p^c](t,s) &= -\sum_{\alpha=1}^{N_r}\left(\frac{2\pi in_{\alpha}c_{n_{\alpha}}}{L_{\alpha}}\sum_{k_1=-\frac{N_1}{2}}^{\frac{N_1}{2}-1}\cdots\sum_{k_{N_r}=-\frac{N_{N_r}}{2}}^{\frac{N_{N_r}}{2}-1}f_{n_1-k_1,\ldots,n_{N_r}-k_{N_r}}[a_{\alpha}](t,s)f_{k_1,\ldots,k_{N_r}}[p^c](t,s)\right) \nonumber \\
        &-\sum_{\alpha=1}^{N_r}\sum_{\substack{\beta=1 \\ \beta\neq\alpha}}^{N_r}\left(\frac{4\pi^2n_{\alpha}n_{\beta}c_{n_{\alpha}}c_{n_{\beta}}}{L_{\alpha}L_{\beta}}\sum_{k_1=-\frac{N_1}{2}}^{\frac{N_1}{2}-1}\cdots\sum_{k_{N_r}=-\frac{N_{N_r}}{2}}^{\frac{N_{N_r}}{2}-1}f_{n_1-k_1,\ldots,n_{N_r}-k_{N_r}}[D_{\alpha\beta}](t,s)f_{k_1,\ldots,k_{N_r}}[p^c](t,s)\right) \nonumber \\
        &-\sum_{\alpha=1}^{N_r}\left(\frac{4\pi^2n_{\alpha}^2}{L_{\alpha}^2}\sum_{k_1=-\frac{N_1}{2}}^{\frac{N_1}{2}-1}\cdots\sum_{k_{N_r}=-\frac{N_{N_r}}{2}}^{\frac{N_{N_r}}{2}-1}f_{n_1-k_1,\ldots,n_{N_r}-k_{N_r}}[D_{\alpha\alpha}](t,s)f_{k_1,\ldots,k_{N_r}}[p^c](t,s)\right)
    \end{align}
\end{figure*}


\section{Coefficient matrices of the ODE systems}

The coefficient matrix $B$ in \eqref{eqn:FPDiscODE} is given at \eqref{eqn:DiscMatrix}, where the variable $e_{i,j,m,n}$ is defined as
\begin{equation}
    e_{i,j,m,n} = \begin{cases}
        \kappa(r_m,s_n,r_i,s_j)\lambda(r_m,s_n)\Delta(r_m) - \lambda(r_i,s_j), \\
            \hspace{1.4in} \text{if} \; m=i \; \text{and} \; n=j \\
        \kappa(r_m,s_n,r_i,s_j)\lambda(r_m,s_n)\Delta(r_m), \\
            \hspace{1.4in} \text{otherwise}.
    \end{cases}
\end{equation}

\begin{figure*}
    \begin{equation} \label{eqn:DiscMatrix}
        B = \left[\begin{matrix}
            e_{1,1,1,1} & \ldots & e_{1,1,N_g,1} &  & e_{1,1,1,N_s} & \ldots & e_{1,1,N_g,N_s} \\
            \vdots & \ddots & \vdots & \qquad \ldots \qquad\qquad & \vdots & \ddots & \vdots \\
            e_{N_g,1,1,1} & \ldots & e_{N_g,1,N_g,1} & & e_{N_g,1,1,N_s} & \ldots & e_{N_g,1,N_g,N_s} \\ \\
            & \vdots & & \ddots & & \vdots & \\ \\
            e_{1,N_s,1,1} & \ldots & e_{1,N_s,N_g,1} & & e_{1,N_s,1,N_s} & \ldots & e_{1,N_s,N_g,N_s} \\
            \vdots & \ddots & \vdots & \qquad \ldots \qquad\qquad & \vdots & \ddots & \vdots \\
            e_{N_g,N_s,1,1} & \ldots & e_{N_g,N_s,N_g,1} & & e_{N_g,N_s,1,N_s} & \ldots & e_{N_g,N_s,N_g,N_s}
        \end{matrix}\right].
    \end{equation}
\end{figure*}

\end{document}